\documentclass[11pt,twoside,reqno]{amsart}
\allowdisplaybreaks
\usepackage{amsmath,amstext,amssymb,epsfig,multicol,enumerate}
\usepackage{stmaryrd}
\usepackage{graphicx}
\usepackage{mathrsfs}
\usepackage{longtable}
\calclayout
\makeatletter
\renewcommand{\@seccntformat}[1]{\bf\csname the#1\endcsname.}
\renewcommand{\section}{\@startsection{section}{1}
	\z@{.7\linespacing\@plus\linespacing}{.5\linespacing}
	{\normalfont\upshape\bfseries\centering}}
\renewcommand{\@biblabel}[1]{\@ifnotempty{#1}{#1.}}
\makeatother
\theoremstyle{plain}
\newtheorem{thm}{Theorem}[section]

\newtheorem{prop}[thm]{Proposition}
\newtheorem{cor}[thm]{Corollary}

\theoremstyle{definition}

\newtheorem{defn}[thm]{Definition}

\newtheorem{rem}{Remark}[section]

\usepackage{cancel}
\usepackage[parfill]{parskip}
\usepackage[german]{varioref}
\usepackage[all]{xy}
\usepackage{color}

\def \>{\succ}
\def \<{\prec}

\begin{document}	
  \title[Ahmed Zahari Abdou\textsuperscript{1}, Kol B\'eatrice GAMOU \textsuperscript{2}, Ibrahima BAKAYOKO\textsuperscript{3} ]{Computational Methods for Rota-type operators on 4-dimensionnal nilpotent Leibniz Algebras.}
	\author{Ahmed Zahari Abdou\textsuperscript{1}, Kol B\'eatrice GAMOU\textsuperscript{2},  Ibrahima BAKAYOKO\textsuperscript{3}}

		\address{\textsuperscript{1}IRIMAS-Department of Mathematics, Faculty of Sciences, University of Haute Alsace, Mulhouse, France}
		\address{\textsuperscript{2}D\'epartement de Math\'ematiques, Universit\'e Gamal Abdel Nasser, Conakry, Guinea.}
        \address{\textsuperscript{2}D\'epartement de Math\'ematiques, Universit\'e de N'Z\'er\'ekor\'e, Guinea.}
	
\email{\textsuperscript{1}abdou-damdji.ahmed-zahari@uha.fr}
\email{\textsuperscript{2}kolbeatrice18@gmail.com}
\email{\textsuperscript{3}ibrahimabakayoko27@gmail.com}
	
	\keywords{Leibniz algebra,  Nilpotent-Leibniz algebra, compatibility, classification.
}
	\subjclass[2010]{16 D70}
	
	\date{\today}
	\begin{abstract}  
A compatible nilpotent Leibniz algebra is a vector space equipped with two multiplication structures that interact in a certain natural way. This article presents the classification of these algebras with dimensions less than four, as well as the classifications of their corresponding Rota-type operators.
\end{abstract}
\maketitle \section{ Introduction}\label{introduction}

Nilpotent Leibniz algebras constitute a natural generalization of Lie algebras, extending some of their fundamental structural properties. As with Lie algebras, a Leibniz algebra is said to be \textit{nilpotent} if its lower central series terminates at zero after a finite number of steps. Since the 2000s, nilpotent Leibniz algebras have attracted growing interest, particularly in the study of their classification, structure, and invariants, especially in the context of low-dimensional cases, filiform algebras, and graded structures.
Leibniz algebras were originally introduced by \textsc{Bloh} in the 1960s under the name \textit{D-algebras}. They were later rediscovered in the 1990s by \textsc{Jean-Louis Loday} in connection with his work on periodicity in algebraic K-theory. Since then, these algebras have received increasing attention due to their broad range of applications in algebraic topology, theoretical physics, and noncommutative geometry.

The study of mathematical structures is often better understood through the operators associated with them. For example, in Galois theory, fields are studied via their automorphisms; in analysis, functions are explored through their derivatives; and in geometry, manifolds are analyzed through their vector fields.

About fifty years ago, several important operators were discovered within the contexts of analysis, probability, and combinatorics. Among these operators are the centroid elements \cite{RM2}, averaging operators, Reynolds operators, Leroux's TD-operator, the Nijenhuis operator, and the Rota-Baxter operator \cite{LG, MD}.

The Rota-Baxter operator finds its origin in the work of G. Baxter \cite{GB, AM} related to Spitzer's identity \cite{F} in fluctuation theory. Rota-Baxter algebras (associative algebras equipped with such an operator) are now used in many areas of mathematics, notably in algebra, number theory, operad theory, and combinatorics \cite{MA, CA1, CA2, PC, GCB}. In mathematical physics, they appear as the operatorial form of the classical Yang-Baxter equation \cite{CA1}, or as the fundamental algebraic structure in the theory of renormalization of quantum fields developed by Connes and Kreimer \cite{CK}.

In the context of non-associative algebras, Rota-Baxter operators are used to construct new operators of the same or different types. For example, the Nijenhuis operator on an associative algebra was introduced in \cite{CJ} in the study of quantum bi-Hamiltonian systems. Regarding Lie algebras, the concept of a Nijenhuis operator originates from the Nijenhuis tensor, introduced in the study of pseudo-complex manifolds, and is related to the well-known notions of the Schouten-Nijenhuis bracket, the Frölicher-Nijenhuis bracket \cite{FN}, and the Nijenhuis-Richardson bracket.

Moreover, the associative analogue of the Nijenhuis relation can be seen as a homogeneous version of the Rota-Baxter relation \cite{PL}.

The (co)homology and homotopy theories of Leibniz algebras have been developed in a number of foundational works \cite{MN}. In recent years, Leibniz algebras have been explored from various perspectives, with applications both in mathematics and physics. These include the integration of Leibniz algebras \cite{T}, deformation quantization \cite{S}, and their role as underlying algebraic structures in embedding tensors, Courant algebroids \cite{T1}, and higher gauge theories 
\cite{S1, SZB, AB, A}.

The classification of algebraic structures and their invariants is an important and ongoing research area in mathematics and physics (see, for example, 
$(L, \llbracket-,-\rrbracket)$ Classical settings include Leibniz and Lie. Specifically, a Leibniz algebra $(L,\llbracket-,-\rrbracket )$ consists of a vector space $L$ equipped with a bilinear map
$\llbracket-,-\rrbracket : L \otimes L\rightarrow L$ that satisfies the associativity condition
\begin{equation*}
\llbracket x,\llbracket y, z\rrbracket\rrbracket=\llbracket\llbracket x, y\rrbracket, z\rrbracket-\llbracket \llbracket x, z\rrbracket, y\rrbracket
\end{equation*}

for all $x, y,z\in L$. In $\llbracket-,-\rrbracket$, the authors classify nilpotent Leibniz algebras over the complex field with dimensions less than four. In a related note, invariant structures such as derivations and centroids have been studied for various classes of algebras in numerous works $(L, \llbracket-,-\rrbracket)$. The algebra of derivations, for example, is useful in algebraic and geometric classification problems.

In the present paper, we are concerned with  Leibniz algebras, which are characterized by a pair of associative-algebraic structures over a common vector space that interact in a nice way. Specifically, a Leibniz algebra is a triple $(L, \llbracket-,-\rrbracket, \left\{-,-\right\})$ for which $(L, \llbracket-,-\rrbracket)$ is the two associative algebras that satisfy the compatibility axiom.

The paper is structured as follows. We first classify Leibniz algebras of dimensions less than four. We then classify automorphisms for all algebras in these dimensions. For a selection of cases, we classify additional invariants such as generalized derivations. The computations for our classifications were done using Mathematica. Throughout, we work on the complex field.

This paper focuses on computational methods for finding Rota-types of 4-dimensional complex Leibniz algebras. We present a systematic approach to compute Rota-Baxter operator, Nenjenhuis operator, Reynolds operator, and Overaging operator using algorithms that can be implemented in symbolic computation software such as Mathematica.

The paper is structured as follows.
\begin{itemize}
    \item Section 2: provides basic definitions and properties of Leibniz algebras, Rota-Baxter operator, Nenjenhuis operator, Reynolds operator and Overaging operator.
    \item Section 3: introduces a step-by-step method for computing derivations of 4-dimensional Leibniz algebras.
    \item The conclusion summarizes the results and suggests directions for future research on compatible Leibniz algebras.
\end{itemize}

This study contributes to a deeper understanding of the computational aspects of Leibniz algebras and offers tools for analyzing their Rota-types. These methods are particularly beneficial for researchers investigating the structure and symmetry of low-dimensional algebraic structures.

 \section{Preliminaries}
We recall here the classification of complex Leibniz algebras of dimension 4, as established in \cite{AOR}.
From these results, we proceed to the classification of complex Leibniz algebras of the same dimension, considered as pairs of compatible Leibniz algebras.
By solving the associated system of equations, we obtain an isomorphism-neutral classification of these algebras using formal calculus methods.
In this subsection, we explore the classification of complex Leibniz algebras in dimension 4, building on the results presented in \cite{AOR}.



\begin{defn}\label{d1}
 Let $L$ be a vector space over a field $\mathbb{K}$ and $x,y,z\in L$. 
 \begin{itemize}
     \item [(i)] A Leibniz algebra structure on $L$ is a bilinear map $\llbracket-,-\rrbracket : L\otimes L\rightarrow L$
     \begin{equation}
\llbracket x,\llbracket y, z\rrbracket\rrbracket=\llbracket\llbracket x, y\rrbracket, z\rrbracket-\llbracket \llbracket x, z\rrbracket, y\rrbracket
\end{equation}
\item [(ii)] A Leibniz algebra $L$ is said to be nilpotent if there exists a natural  $s\in \mathbb{N}$, such that $l^s=0.$
\item [(iii)] Two Leibniz algebras $(L, \llbracket-,-\rrbracket_1)$ and $(L, \llbracket-,-\rrbracket_2)$ are called compatible if for any 
$\lambda_1,\lambda_2\in \mathbb{K},$ the following bracket 
 \begin{equation}\label{d}
\llbracket x,y\rrbracket=\lambda_1\llbracket x, y\rrbracket_1+\lambda_2\llbracket x, y\rrbracket_2,\quad \forall \, x,y\in L,
\end{equation}
define a Leibniz algebra structure on L.
 \end{itemize}
\end{defn}

\begin{rem}
The bracket (\ref{d}) defines a Leibniz algebra structure on $L$ if and only if 
\begin{equation}\label{dc}
	\begin{aligned}
	\llbracket x,\llbracket y, z\rrbracket_1\rrbracket_2+\llbracket x,\llbracket y, z\rrbracket_2\rrbracket_1=&\llbracket\llbracket x, y\rrbracket_1, z\rrbracket_2+\llbracket\llbracket x, y\rrbracket_2, z\rrbracket_1-\llbracket\llbracket x,z\rrbracket_1,y\rrbracket_2-\llbracket\llbracket x,z\rrbracket_2,y\rrbracket_1
	\end{aligned}
	\end{equation}
for $x, y, z\in L$.
\end{rem}

\begin{defn}\label{dcl}
The triple $ (L, \llbracket-,-\rrbracket_1, \llbracket-,-\rrbracket_2)$ is said to be a 
compatible Leibniz algebra if $(\llbracket-,-\rrbracket_1)$ and $(\llbracket-,-\rrbracket_2)$
are both Leibniz algebras and (\ref{dc}) holds.
\end{defn}

\begin{defn}
A homorphism $\theta : (L, \llbracket-,-\rrbracket)\rightarrow  (L', \llbracket-,-\rrbracket')$ of Leibniz algebras is $K$-linear maps satisfying $\theta(\llbracket x, y\rrbracket)=\llbracket \theta(x), \theta(y)\rrbracket',\, \forall\, x, y\in L.$
\end{defn}

\begin{defn}\label{d3}
Let $(L, \llbracket-,-\rrbracket)$ be a Leibniz algebra. A linear map $r : L\rightarrow L$ is said to be a Rota-Baxter operator  if 
\begin{align}
   \llbracket r(x), r(y)\rrbracket &= r\Bigg(\llbracket r(x), y \rrbracket +  \llbracket x, r(y) \rrbracket + \lambda  \llbracket x, y  \rrbracket\Bigg)\label{eq3} , 
\end{align}
for any $\, x, y\in L.$
\end{defn}

\begin{defn}\label{d4}
Let $(L, \llbracket-,-\rrbracket)$ be a Leibniz algebra. A linear map $\kappa : L\rightarrow L$ is said to be a Nenjenhuis operator  if 
\begin{align}
  \llbracket \kappa(x), \kappa(y)\rrbracket&=\kappa\Bigg(\llbracket \kappa(x), y \rrbracket +  \llbracket x, \kappa(y) \rrbracket -\kappa( \llbracket x,y \rrbracket)\Bigg)\label{eq4},
\end{align}
for any $\, x, y\in L.$
\end{defn}

\begin{defn}\label{d5}
Let $(L, \llbracket-,-\rrbracket)$ be a Leibniz algebra. A linear map $\alpha : L\rightarrow L$ is said to be a Reynold operator  if 
\begin{align}
   \llbracket \alpha(x), \alpha(y)\rrbracket &= \alpha\Bigg(\llbracket x, \alpha(y) \rrbracket +  \llbracket \alpha(x), y \rrbracket -  \llbracket \alpha(x), \alpha(y) \rrbracket\Bigg)\label{eq5}, 
\end{align}
for any $\, x, y\in L.$
\end{defn}

\begin{defn}\label{d6}
Let $(L, \llbracket-,-\rrbracket)$ be a Leibniz algebra. A linear map $\beta : L\rightarrow L$ is said to be a Overaging operator  if 
\begin{align}
  \beta(\llbracket x, \beta(y)\rrbracket) &= \llbracket \beta(x), \beta(y)\rrbracket=\beta(\llbracket x, \beta(y)\rrbracket)\label{eq6},
\end{align}
for any $\, x, y\in L.$
\end{defn}

\begin{thm}
The isomorphism class of four-dimensional complex nilpotent Leibniz algebras given by  the following representatives.
\begin{itemize}
\item
$L_1$ :
 $\begin{array}{ll}  
\llbracket e_1,e_1\rrbracket=e_2,\quad \llbracket e_2,e_1\rrbracket=e_3,\quad\llbracket e_3,e_1\rrbracket=e_4;
\end{array}$
\item
$L_2$ :
 $\begin{array}{ll}  
\llbracket e_1,e_1\rrbracket=e_3,\quad\rrbracket e_1,e_2\rrbracket=e_4,\quad\rrbracket e_2,e_1\rrbracket=e_3,\quad \rrbracket e_3,e_1\rrbracket=e_4;
\end{array}$
\item
$L_{3}$ :
 $\begin{array}{ll}  
\llbracket e_1,e_1\rrbracket=e_3,\quad\rrbracket e_2,e_1\rrbracket e_3,\quad\rrbracket e_3,e_1\rrbracket=e_4;
\end{array}$
\item
$L_{4}(\mu)$:
 $\begin{array}{ll}  
\llbracket e_1,e_1\rrbracket= e_3,\quad\llbracket e_1,e_2\rrbracket=\mu e_4,\quad\llbracket e_2,e_1\llbracket=e_3,\quad\llbracket e_3,e_2\rrbracket=e_4,\quad\left[e_2,e_2\right]=e_4,\quad \mu\in\left\{0,1\right\};
\end{array}$
\item
$L_{5}$ :
 $\begin{array}{ll}  
\llbracket e_1,e_1\rrbracket=e_3,\quad \rrbracket e_1,e_2\rrbracket=e_4,\quad\rrbracket e_3,e_1\rrbracket=e_4;
\end{array}$
\item
$L_{6}$ :
 $\begin{array}{ll}  
\llbracket e_1,e_1\rrbracket=e_3,\quad \llbracket e_2,e_2\rrbracket=e_4,\quad\llbracket e_3,e_1\rrbracket =e_4;
\end{array}$
\item
$L_{7}$ :
 $\begin{array}{ll}  
\llbracket e_1,e_1\rrbracket=e_4,\quad \llbracket e_1,e_2\rrbracket=-e_3,\quad\llbracket e_1,e_3\rrbracket=-e_4,\quad
\llbracket e_2,e_1\rrbracket=e_3,\quad\llbracket e_3,e_1\rrbracket =e_4 ;
\end{array}$
\item
$L_{8}$ :
 $\begin{array}{ll}  
\llbracket e_1,e_1\rrbracket=e_4,\quad \llbracket e_1,e_2\rrbracket=-e_3+e_4,\quad\llbracket e_1,e_3\rrbracket=-e_4,\quad
\llbracket e_2,e_1\rrbracket=e_3,\quad\llbracket e_3,e_1\rrbracket=e_4 ;
\end{array}$
\item
$L_{9}$ :
 $\begin{array}{ll}  
\llbracket e_1,e_1\rrbracket=e_4,\,\llbracket e_1,e_2\rrbracket=-e_3+2e_4,\,\llbracket e_1,e_3\rrbracket=-e_4,\,
\llbracket e_2,e_1\rrbracket=e_3,\,\llbracket e_2,e_2\rrbracket=e_4,\,\llbracket e_3,e_1\rrbracket=e_4;
\end{array}$
\item
$L_{10}$ :
 $\begin{array}{ll}  
\llbracket e_1,e_1\rrbracket=e_4,\,\llbracket e_1,e_2\rrbracket=-e_3,\,\llbracket e_1,e_3\rrbracket=-e_4,\,
\llbracket e_2,e_1\rrbracket=e_3,\,\llbracket e_2,e_2\rrbracket=e_4,\,\llbracket e_3,e_1\rrbracket=e_4;
\end{array}$
\item
$L_{11}$ :
 $\begin{array}{ll}  
\llbracket e_1,e_1\rrbracket=e_4,\,\llbracket e_1,e_2\rrbracket=e_3,\,\llbracket e_2,e_1\rrbracket=-e_3,\,
\llbracket e_2,e_2\rrbracket=-2e_3+e_4;
\end{array}$
\item
$L_{12}$ :
 $\begin{array}{ll}  
\llbracket e_1,e_2\rrbracket= e_3,\quad\llbracket e_2,e_1\rrbracket=e_4,\quad\llbracket e_2,e_2\rrbracket=-e_3;
\end{array}$
\item
$L_{13}(\mu)$ :
 $\begin{array}{ll}  
\llbracket e_1,e_1\rrbracket=e_3,\,\llbracket e_1,e_2\rrbracket=e_4,\,\llbracket e_2,e_1\rrbracket=-\mu e_3,\,\llbracket e_2,e_2\rrbracket=-e_4,\quad \mu \in \mathbb{C};
\end{array}$
\item
$L_{14}(\mu)$ :
 $\begin{array}{ll}  
\llbracket e_1,e_1\rrbracket=e_4,\,\llbracket e_1,e_2\rrbracket=\mu e_4,\,\llbracket e_2,e_1 \rrbracket=-\alpha e_4,\,\llbracket e_2,e_2\rrbracket=e_4,\,\llbracket e_3,e_3\rrbracket=e_4,\quad \mu \in \mathbb{C};
\end{array}$
\item
$L_{15}$ :
 $\begin{array}{ll}  
\llbracket e_1,e_2\rrbracket=e_4,\,\llbracket e_1,e_3\rrbracket=e_4,\,\llbracket e_2,e_1\rrbracket=-e_4,\,\llbracket e_2,e_2\rrbracket=e_4,\,\llbracket e_3,e_1\rrbracket=e_4;
\end{array}$
\item
$L_{16}$ :
 $\begin{array}{ll}  
\llbracket e_1,e_1\rrbracket=e_4,\,\llbracket e_1,e_2\rrbracket=e_4,\,\llbracket e_2,e_1\rrbracket=-e_4,\,\llbracket e_3,e_3\rrbracket=e_4;
\end{array}$
\item
$L_{17}$ :
 $\begin{array}{ll}  
\llbracket e_1,e_2\rrbracket =e_3,\,\llbracket e_2,e_1\rrbracket =e_4;
\end{array}$
\item
$L_{18}$ :
 $\begin{array}{ll}  
\llbracket e_1,e_2\rrbracket=e_3,\,\llbracket e_2,e_1\rrbracket=-e_3,\,\llbracket e_2,e_2\rrbracket =e_4;
\end{array}$
\item
$L_{19}$ :
 $\begin{array}{ll}  
\llbracket e_2,e_1\rrbracket =e_4,\quad\llbracket e_2,e_2\rrbracket= e_3;
\end{array}$
\item
$L_{20}(\mu)$ :
 $\begin{array}{ll}  
\llbracket e_1,e_2\rrbracket =e_4,\,\llbracket e_2,e_1\rrbracket =\frac{1+\mu}{1-\mu}e_4,\,\llbracket e_2,e_2\rrbracket =e_3,\mu\in\mathbb{C}\setminus\{1\};
\end{array}$
\item
$L_{21}$ :
 $\begin{array}{ll}  
\llbracket e_1,e_2\rrbracket=e_4,\quad\llbracket e_2,e_1\rrbracket=-e_4,\quad \llbracket e_3,e_3\rrbracket=e_4.
\end{array}$
\end{itemize}
\end{thm}

\newpage
\section{Rota-type operators}
\subsection{Rota-Baxter operator}\,

 Let $(L, \llbracket-,-\rrbracket)$   be an $n-$dimensional Leibniz algebra with basis $\{e_i\}\, (1 \leq i \leq n)$ and let $r$ be a Rota-Baxter operator on $L$
For any $i,j,k\in\mathbb{N}$, $1 \leq i,j,k \leq n$, let us put
$$\llbracket e_i,e_j\rrbracket=\sum_{k=1}^nc_{ij}^ke_k\quad \text{and}\quad r(e_i)=\sum_{j=1}^nr_{ji}e_j.$$
Then, in term of basis elements, equation (\ref{eq3}) is equivalent to
\begin{align*}
   \sum^{n}_{k=1}\sum^{n}_{p=1}r_{ki}r_{pj} C_{kp}^{p}= 
   \sum^{n}_{k=1} \sum^{n}_{p=1}\Bigg(r_{ki}C_{kj}^{p}r_{qp}+r_{kj}C_{ik}^{p}r_{qp}\Bigg)+\lambda\sum_{p=1}^nC_{ij}^pr_{qp},
\end{align*} for $i,j,q=1,2,\dots,n$.

\begin{prop}
The description of the Rota-Baxter operator of every 4-dimensional  Leibniz algebra is given below.
\begin{itemize}
    \item[$L_1$ : ]
    $\left(\begin{array}{ccccc}
0&0&0&0\\
0&0&0&0\\
r_{31}&r_{32}&r_{33}&0\\
r_{41}&r_{42}&r_{43}&0\\
\end{array}
\right),\quad
\left(\begin{array}{ccccc}
0&0&0&0\\
r_{21}&r_{22}&0&0\\
r_{31}&r_{32}&0&0\\
r_{41}&r_{42}&0&0\\
\end{array}
\right)
,\quad
\left(\begin{array}{ccccc}
0&0&0&0\\
0&r_{22}&0&0\\
0&r_{32}&0&0\\
r_{41}&r_{42}&r_{43}&-\frac{r_{22}r_{43}}{r_{32}}\\
\end{array}
\right),\\
\left(\begin{array}{ccccc}
0&0&0&0\\
r_{21}&r_{22}&0&0\\
0&0&0&0\\
r_{41}&r_{42}&0&r_{44}
\end{array}
\right),\quad
\left(\begin{array}{ccccc}
0&0&0&0\\
r_{21}&0&0\\
0&0&0&0\\
r_{41}&r_{42}&0&r_{44}\\
\end{array}
\right)
,\quad
\left(\begin{array}{ccccc}
0&0&0&0\\
0&0&0&0\\
0&0&0&0\\
r_{41}&r_{42}&r_{43}&r_{44}\\
\end{array}
\right)
;
$

\item[$L_2$ : ]
    $\left(\begin{array}{ccccc}
0&0&0&0\\
0&0&0&0\\
r_{31}&r_{32}&r_{33}&0\\
r_{41}&r_{42}&r_{43}&0\\
\end{array}
\right),\quad
\left(\begin{array}{ccccc}
0&0&0&0\\
0&0&0&0\\
r_{31}&r_{32}&0&0\\
r_{41}&r_{42}&r_{43}&0\\
\end{array}
\right)
,\quad
\left(\begin{array}{ccccc}
0&0&0&0\\
r_{21}&0&0&0\\
0&0&0&0\\
r_{41}&r_{42}&-r_{44}&r_{44}\\
\end{array}
\right),\\
\left(\begin{array}{ccccc}
r_{11}&0&0&0\\
-r_{11}&0&0&0\\
r_{11}+r_{43}&r_{43}&0&0\\
r_{41}&r_{42}&r_{43}&0
\end{array}
\right),\quad
\left(\begin{array}{ccccc}
0&0&0&0\\
0&0&0\\
0&0&0&0\\
r_{41}&r_{42}&0&r_{44}\\
\end{array}
\right)
;
$

\item[$L_3$ : ]
    $\left(\begin{array}{ccccc}
0&0&0&0\\
0&0&0&0\\
r_{31}&r_{32}&r_{33}&0\\
r_{41}&r_{42}&r_{43}&0\\
\end{array}
\right),\quad
\left(\begin{array}{ccccc}
0&0&0&0\\
r_{21}&r_{22}&0&0\\
0&r_{32}&0&0\\
r_{41}&r_{42}&0&0\\
\end{array}
\right)
,\quad
\left(\begin{array}{ccccc}
0&0&0&0\\
r_{21}&r_{22}&0&0\\
0&0&0&0\\
r_{41}&r_{42}&0&r_{44}\\
\end{array}
\right),\\
\left(\begin{array}{ccccc}
r_{11}&r_{12}&0&0\\
-r_{11}&-r_{12}&0&0\\
r_{31}&r_{33}&0&0\\
r_{41}&r_{42}&r_{43}&0
\end{array}
\right),\quad
\left(\begin{array}{ccccc}
0&0&0&0\\
0&0&0\\
0&0&0&0\\
r_{41}&r_{42}&r_{43}&r_{44}\\
\end{array}
\right)
,\quad
\left(\begin{array}{ccccc}
0&0&0&0\\
r_{21}&r_{22}&0&r_{24}\\
0&0&0&0\\
r_{41}&r_{42}&r_{43}&r_{44}\\
\end{array}
\right)
;
$

\item[$L_4$ : ]
    $\left(\begin{array}{ccccc}
0&0&0&0\\
r_{21}&0&0&0\\
r_{31}&r_{32}&0&0\\
r_{41}&r_{42}&r_{43}&0\\
\end{array}
\right),\quad
\left(\begin{array}{ccccc}
0&0&0&0\\
0&0&0&0\\
r_{31}&r_{32}&r_{33}&0\\
r_{41}&r_{42}&r_{43}&0\\
\end{array}
\right)
,\quad
\left(\begin{array}{ccccc}
0&0&0&0\\
0&0&0&0\\
0&0&0&0\\
r_{41}&r_{42}&r_{43}&r_{44}\\
\end{array}
\right)
;
$

\item[$L_5$ : ]
    $\left(\begin{array}{ccccc}
0&0&0&0\\
r_{21}&r_{22}&r_{23}&0\\
r_{31}&r_{32}&r_{33}&0\\
r_{41}&r_{42}&r_{43}&0\\
\end{array}
\right),\quad
\left(\begin{array}{ccccc}
0&0&0&0\\
r_{21}&r_{22}&-r_{33}&0\\
r_{31}&r_{32}&r_{33}&0\\
r_{41}&r_{42}&r_{43}&0\\
\end{array}
\right)
,\quad
\left(\begin{array}{ccccc}
0&0&0&0\\
r_{21}&r_{22}&0&0\\
r_{31}&r_{32}&r_{33}&0\\
r_{41}&r_{42}&r_{43}&0\\
\end{array}
\right),
$

$
\quad
\left(\begin{array}{ccccc}
0&0&0&0\\
r_{21}&0&0&0\\
-r_{21}&0&0&0\\
r_{41}&r_{42}&r_{43}&r_{44}\\
\end{array}
\right),
\,
\left(\begin{array}{ccccc}
r_{11}&0&0&0\\
r_{21}&\frac{1}{2}r_{11}&0&0\\
r_{31}&0&\frac{1}{2}r_{11}&0\\
r_{41}&r_{42}&r_{43}&\frac{1}{2}(r_{21}+r_{31})\\
\end{array}
\right)
;
$

\item[$L_6$ : ]
    $\left(\begin{array}{ccccc}
0&0&0&0\\
0&0&0&0\\
r_{31}&r_{32}&r_{33}&0\\
r_{41}&r_{42}&r_{43}&0\\
\end{array}
\right),\quad
\left(\begin{array}{ccccc}
r_{11}&0&0&0\\
0&0&0&0\\
r_{31}&0&\frac{1}{2}r_{32}&0\\
r_{41}&r_{42}&\frac{1}{2}r_{31}&\frac{1}{3}r_{11}&\\
\end{array}
\right)
,\quad
\left(\begin{array}{ccccc}
r_{11}&0&0&0\\
0&\frac{2}{3}r_{11}&0&0\\
r_{31}&0&\frac{1}{3}r_{11}&0\\
r_{41}&r_{42}&\frac{1}{3}r_{31}&\frac{1}{3}r_{11}\\
\end{array}
\right),
$

$
\quad
\left(\begin{array}{ccccc}
0&0&0&0\\
0&0&0&0\\
0&0&0&0\\
r_{41}&r_{42}&r_{43}&r_{44}\\
\end{array}
\right),
\,
\left(\begin{array}{ccccc}
0&0&0&0\\
0&2r_{44}&0&0\\
0&0&0&0\\
r_{41}&r_{42}&r_{43}&r_{44}\\
\end{array}
\right)
;
$

\item[$L_7$ : ]
    $\left(\begin{array}{ccccc}
0&0&0&0\\
r_{21}&0&0&0\\
r_{31}&r_{32}&r_{33}&0\\
r_{41}&r_{42}&r_{43}&0\\
\end{array}
\right),\quad
\left(\begin{array}{ccccc}
0&0&0&0\\
r_{21}&r_{22}&0&0\\
r_{31}&r_{32}&0&0\\
r_{41}&r_{42}&-\frac{r_{32}r_{44}}{r_{22}}&r_{44}\\
\end{array}
\right)
,\quad
\left(\begin{array}{ccccc}
0&0&0&0\\
r_{21}&r_{22}&0&r_{24}\\
r_{31}&0&0&0\\
r_{41}&r_{42}&0&r_{44}\\
\end{array}
\right),
$

$
\quad
\left(\begin{array}{ccccc}
0&0&0&0\\
r_{21}&0&0&0\\
r_{31}&r_{32}&0&0\\
r_{41}&r_{42}&r_{43}&0\\
\end{array}
\right),
\quad
\left(\begin{array}{ccccc}
0&0&0&0\\
r_{21}&0&0&0\\
r_{31}&0&0&0\\
r_{41}&r_{42}&r_{43}&r_{44}\\
\end{array}
\right)
,
\quad
\left(\begin{array}{ccccc}
0&0&0&0\\
r_{21}&r_{22}&0&0\\
r_{31}&r_{32}&0&0\\
r_{41}&r_{42}&0&0\\
\end{array}
\right)
;
$

\item[$L_8$ : ]
    $\left(\begin{array}{ccccc}
0&0&0&0\\
r_{21}&0&0&0\\
r_{31}&r_{32}&r_{33}&0\\
r_{41}&r_{42}&r_{43}&0\\
\end{array}
\right),\quad
\left(\begin{array}{ccccc}
0&0&0&0\\
0&0&0&0\\
r_{31}r_{32}&r_{33}&0\\
r_{41}&r_{42}&r_{43}&0\\
\end{array}
\right)
,\quad
\left(\begin{array}{ccccc}
0&r_{12}&0&0\\
0&-r_{12}&0&0\\
r_{43}&r_{32}&0&0\\
r_{41}&r_{42}&r_{43}&0\\
\end{array}
\right),
$

$
\quad
\left(\begin{array}{ccccc}
r_{11}&0&0&0\\
-r_{11}&0&0&0\\
r_{31}&r_{32}&0&0\\
r_{41}&r_{42}&r_{43}&0\\
\end{array}
\right),
\qquad
\left(\begin{array}{ccccc}
0&0&0&0\\
r_{21}&r_{22}&0&0\\
r_{31}&r_{32}&0&0\\
r_{41}&r_{42}&0&0\\
\end{array}
\right)
;
$

\item[$L_9$ : ]
    $\left(\begin{array}{ccccc}
0&0&0&0\\
0&0&0&0\\
r_{31}&r_{32}&r_{33}&0\\
r_{41}&r_{42}&r_{43}&0\\
\end{array}
\right),\quad
\left(\begin{array}{ccccc}
0&0&0&0\\
0&0&0&0\\
r_{31}&0&0&0\\
r_{41}&r_{42}&r_{43}&r_{44}\\
\end{array}
\right)
,\quad
\left(\begin{array}{ccccc}
r_{11}&r_{11}&0&0\\
-r_{11}&-r_{11}&0&0\\
r_{31}&r_{32}&0&0\\
r_{41}&r_{42}&r_{43}&0\\
\end{array}
\right),
$

$
\quad
\left(\begin{array}{ccccc}
0&0&0&0\\
2r_{44}&2r_{43}&0&0\\
r_{31}&r_{32}&0&0\\
r_{41}&r_{42}&\frac{1}{2}r_{32}+r_{44}&r_{44}\\
\end{array}
\right),
\,
\left(\begin{array}{ccccc}
r_{11}&0&0&0\\
-\frac{1}{2}r_{11}&\frac{1}{2}r_{11}&0&0\\
r_{31}&r_{32}&\frac{1}{3}r_{11}&0\\
r_{41}&r_{42}&\frac{1}{12}(-r_{11}+6r_{32})&\frac{1}{4}r_{11}\\
\end{array}
\right)
;
$

\item[$L_{10}$ : ]
    $\left(\begin{array}{ccccc}
0&0&0&0\\
0&0&0&0\\
r_{31}&r_{32}&r_{33}&0\\
r_{41}&r_{42}&r_{43}&0\\
\end{array}
\right),\quad
\left(\begin{array}{ccccc}
0&0&0&0\\
r_{21}&0&0&0\\
\frac{r_{21}}{2r_{44}}&-r_{21}&0&0\\
r_{41}&r_{42}&r_{43}&r_{44}\\
\end{array}
\right)
,\quad
\left(\begin{array}{ccccc}
0&0&0&0\\
0&0&0&0\\
0&0&0&0\\
r_{41}&r_{42}&r_{43}&r_{44}\\
\end{array}
\right),
$

$
\quad
\left(\begin{array}{ccccc}
0&r_{12}&0&0\\
0&0&0&0\\
0&0&0&0\\
r_{41}&r_{42}&r_{43}&0\\
\end{array}
\right),
\quad
\left(\begin{array}{ccccc}
r_{11}&0&0&0\\
-r_{11}&0&0&0\\
-r_{11}&0&0&0\\
r_{41}&r_{42}&0&0\
\end{array}
\right)
,
\quad
\left(\begin{array}{ccccc}
r_{11}&\frac{1}{2}r_{11}&0&0\\
r_{21}&0&0&0\\
\frac{-r_{11}^2-2r_{21}^2}{3r_{11}}&-\frac{1}{2}r_{21}&\frac{1}{3}r_{11}&0\\
r_{41}&r_{42}&0&\frac{1}{4}r_{11}
\end{array}
\right)
;
$

\item[$L_{11}$ : ]
    $\left(\begin{array}{ccccc}
0&0&0&0\\
0&0&0&0\\
r_{31}&r_{32}&r_{33}&r_{34}\\
r_{41}&r_{42}&r_{43}&r_{44}\\
\end{array}
\right),\quad
\left(\begin{array}{ccccc}
0&0&0&0\\
0&0&0&0\\
r_{31}&r_{32}&r_{33}&3r_{33}\\
r_{41}&r_{42}&r_{43}&r_{44}\\
\end{array}
\right)
,\quad
\left(\begin{array}{ccccc}
\frac{r_{21}^2}{r_{22}}&r_{21}&0&0\\
r_{21}&r_{22}&0&0\\
r_{31}&r_{32}&0&-r_{22}\\
r_{41}&r_{42}&0&\frac{r_{21}^2+r_{22}^2}{2r_{22}}\\\
\end{array}
\right),
$

$
\quad
\left(\begin{array}{ccccc}
r_{22}&0&0&0\\
0&r_{22}&0&0\\
r_{31}&r_{32}&\frac{1}{2}r_{22}&0\\
r_{41}&r_{42}&0&\frac{1}{2}r_{22}\\
\end{array}
\right),
\quad
\left(\begin{array}{ccccc}
r_{11}&0&0&0\\
0&0&0&0\\
r_{31}&r_{32}&0&0\\
r_{41}&r_{42}&0&\frac{1}{2}r_{11}\
\end{array}
\right)
,
\quad
\left(\begin{array}{ccccc}
0&0&0&0\\
0&r_{22}&0&0\\
r_{41}&r_{32}&0&-r_{22}\\
r_{41}&r_{42}&0&\frac{1}{2}r_{11}
\end{array}
\right)
;
$

\item[$L_{12}$ : ]
    $\left(\begin{array}{ccccc}
r_{11}&0&0&0\\
0&0&0&0\\
r_{31}&r_{32}&0&0\\
r_{41}&r_{42}&0&0\\
\end{array}
\right),\quad
\left(\begin{array}{ccccc}
r_{22}&0&0&0\\
0&r_{22}&0&0\\
r_{31}&r_{32}&\frac{1}{2}r_{22}&0\\
r_{41}&r_{42}&0&\frac{1}{2}r_{22}\\
\end{array}
\right)
,\quad
\left(\begin{array}{ccccc}
0&r_{12}&0&0\\
0&0&0&0\\
r_{31}&r_{32}&r_{33}&-r_{33}\\
r_{41}&r_{42}&0&0\\\
\end{array}
\right),
$

$
\quad
\left(\begin{array}{ccccc}
0&0&0&0\\
0&0&0&0\\
r_{31}&r_{32}&r_{22}&r_{34}\\
r_{41}&r_{42}&r_{34}&r_{44}\\
\end{array}
\right),
\quad
\left(\begin{array}{ccccc}
0&r_{12}&0&0\\
0&0&0&0\\
r_{31}&r_{32}&r_{22}&-r_{33}\\
r_{41}&r_{42}&r_{34}&-r_{34}\\
\end{array}
\right)
;
$

\item[$L_{13}$ : ]
    $\left(\begin{array}{ccccc}
0&0&0&0\\
0&0&0&0\\
r_{31}&r_{32}&r_{22}&r_{34}\\
r_{41}&r_{42}&r_{34}&r_{44}\\
\end{array}
\right),\quad
\left(\begin{array}{ccccc}
0&0&0&0\\
0&0&0&0\\
r_{31}&r_{32}&r_{22}&0\\
r_{41}&r_{42}&r_{34}&r_{44}\\
\end{array}
\right)
,\quad
\left(\begin{array}{ccccc}
0&0&0&0\\
0&0&0&0\\
r_{31}&r_{32}&r_{22}&0\\
r_{41}&r_{42}&r_{34}&0\\
\end{array}
\right)
;
$

\item[$L_{14}$ : ]
    $\left(\begin{array}{ccccc}
0&0&0&0\\
0&0&0&0\\
0&0&0\\
r_{41}&r_{42}&r_{34}&r_{44}\\
\end{array}
\right),\quad
\left(\begin{array}{ccccc}
0&0&0&0\\
0&0&r_{23}&0\\
r_{31}&0&-ir_{23}&0\\
r_{41}&r_{42}&r_{34}&0\\
\end{array}
\right)
,\quad
\left(\begin{array}{ccccc}
0&0&0&0\\
0&0&0&0\\
0&0&0&0\\
r_{41}&r_{42}&r_{34}&0\\
\end{array}
\right)
;
$

 $\left(\begin{array}{ccccc}
0&r_{12}&r_{13}&0\\
0&0&0&0\\
0&ir_{12}&ir_{13}&0\\
r_{41}&r_{42}&r_{34}&0\\
\end{array}
\right),\quad
\left(\begin{array}{ccccc}
0&0&0&0\\
0&r_{12}&r_{22}&0\\
ir_{21}&ir_{22}&0&0\\
r_{41}&r_{42}&r_{34}&0\\
\end{array}
\right)
;
$

\item[$L_{15}$ : ]
    $\left(\begin{array}{ccccc}
0&0&0&0\\
0&0&0&0\\
0&0&0\\
r_{41}&r_{42}&r_{34}&r_{44}\\
\end{array}
\right),\quad
\left(\begin{array}{ccccc}
0&0&0&0\\
0&0&0&0\\
r_{31}&0&-ir_{23}&0\\
r_{41}&r_{42}&r_{34}&0\\
\end{array}
\right)
,\quad
\left(\begin{array}{ccccc}
r_{11}&r_{12}&r_{13}&0\\
0&0&0&0\\
0&0&0&0\\
r_{41}&r_{42}&r_{34}&0\\
\end{array}
\right)
;
$

 $\left(\begin{array}{ccccc}
0&0&0&0\\
r_{21}&0&0&0\\
\frac{r_{21}^2}{2r_{44}}&-r_{21}&0&0\\
r_{41}&r_{42}&r_{34}&r_{44}\\
\end{array}
\right),
\left(\begin{array}{ccccc}
2r_{44}&0&0&0\\
r_{21}&2r_{44}&0&0\\
-\frac{r_{21}^2}{2r_{44}}&-r_{21}&2r_{22}&0\\
r_{41}&r_{42}&r_{34}&r_{44}\\
\end{array}
\right)
,
\left(\begin{array}{ccccc}
0&0&0&0\\
r_{21}&0&0&0\\
-\frac{r_{21}^2}{2r_{44}}&-r_{21}&0&0\\
r_{41}&r_{42}&r_{34}&r_{44}\\
\end{array}
\right)
;
$

\item[$L_{16}$ : ]
    $\left(\begin{array}{ccccc}
0&0&0&0\\
r_{21}&r_{22}&r_{23}&0\\
0&0&0\\
r_{41}&r_{42}&r_{34}&0\\
\end{array}
\right),
\left(\begin{array}{ccccc}
0&0&0&0\\
r_{21}&0&-r_{44}&0\\
r_{44}&0&2r_{44}&0\\
r_{41}&r_{42}&r_{34}&r_{44}\\
\end{array}
\right)
,
\left(\begin{array}{ccccc}
2r_{44}&0&0&0\\
r_{21}&2r_{44}&0&0\\
0&0&2r_{22}&0\\
r_{41}&r_{42}&r_{34}&r_{44}\\
\end{array}
\right)
,
 \left(\begin{array}{ccccc}
2r_{44}&0&0&0\\
r_{21}&2r_{44}&r_{44}&0\\
-r_{44}&0&0&0\\
r_{41}&r_{42}&r_{34}&r_{44}\\
\end{array}
\right)
;
$

\item[$L_{17}$ : ]
    $\left(\begin{array}{ccccc}
0&0&0&0\\
0&0&0&0\\
r_{31}&r_{32}&r_{33}&r_{34}\\
r_{41}&r_{42}&r_{34}&r_{44}\\
\end{array}
\right),\quad
\left(\begin{array}{ccccc}
0&0&0&0\\
r_{21}&0&0&0\\
r_{31}&r_{32}&-r_{34}&r_{34}\\
r_{41}&r_{42}&r_{34}&-r_{43}\\
\end{array}
\right)
,\quad
\left(\begin{array}{ccccc}
0&r_{12}&0&0\\
0&0&0&0\\
r_{31}&r_{32}&-r_{34}&r_{34}\\
r_{41}&r_{42}&r_{34}&-r_{43}\\
\end{array}
\right)
,
$

 $\left(\begin{array}{ccccc}
0&0&0&0\\
r_{21}&r_{22}&0&0\\
r_{31}&r_{32}&0&0\\
r_{41}&r_{42}&0&0\\
\end{array}
\right)
,\quad
\left(\begin{array}{ccccc}
r_{11}&0&0&0\\
0&r_{22}&0&0\\
r_{31}&r_{32}&\frac{r_{11}r_{22}}{r_{11}+r_{22}}&0\\
r_{41}&r_{42}&0&\frac{r_{11}r_{22}}{r_{11}+r_{22}}
\end{array}
\right)
;
$

\item[$L_{18}$ : ]
    $\left(\begin{array}{ccccc}
r_{11}&r_{12}&0&0\\
0&r_{22}&0&0\\
r_{31}&r_{32}&\frac{r_{11}r_{22}}{r_{11}+r_{22}}&0\\
r_{41}&r_{42}&0&\frac{1}{2}r_{22}\\
\end{array}
\right),\quad
\left(\begin{array}{ccccc}
0&r_{12}&0&0\\
0&0&0&0\\
r_{31}&r_{32}&r_{33}&0\\
r_{41}&r_{42}&r_{43}&0\\
\end{array}
\right)
,\quad
\left(\begin{array}{ccccc}
r_{11}&r_{12}&0&0\\
0&0&0&0\\
r_{31}&r_{32}&0&r_{34}\\
r_{41}&r_{42}&0&r_{44}\\
\end{array}
\right)
,
$

 $\left(\begin{array}{ccccc}
r_{11}&r_{12}&0&r_{14}\\
0&0&0&0\\
r_{31}&r_{32}&0&r_{34}\\
r_{41}&r_{42}&0&r_{44}\\
\end{array}
\right)
,\quad
\left(\begin{array}{ccccc}
r_{11}&r_{12}&0&0\\
0&0&0&0\\
r_{31}&r_{32}&r_{33}&r_{34}\\
r_{41}&r_{42}&r_{43}&r_{44}\\
\end{array}
\right)
;
$

\item[$L_{19}$ : ]
    $\left(\begin{array}{ccccc}
r_{11}&0&0&0\\
0&r_{22}&0&0\\
r_{31}&r_{32}&\frac{1}{2}r_{22}&0\\
r_{41}&r_{42}&0& \frac{r_{11}r_{22}}{r_{11}+r_{22}} 
\end{array}
\right),
\left(\begin{array}{ccccc}
r_{11}&r_{12}&0&0\\
0&r_{22}&0&0\\
r_{31}&r_{32}&\frac{1}{2}r_{22}&0\\
r_{41}&r_{42}&\frac{r_{11}r_{22}}{2(r_{11}+r_{22})}&\frac{r_{11}r_{22}}{r_{11}+r_{22}} 
\end{array}
\right)
,
\left(\begin{array}{ccccc}
r_{11}&r_{12}&r_{13}&0\\
0&0&0&0\\
r_{31}&r_{32}&r_{33}&0\\
r_{41}&r_{42}&r_{43}&0\\
\end{array}
\right)
,
$

 $\left(\begin{array}{ccccc}
 0&r_{12}&0&0\\
0&r_{22}&0&0\\
r_{31}&r_{32}&\frac{1}{2}r_{22}&0\\
r_{41}&r_{42}&\frac{1}{2}r_{22}&0\\
\end{array}
\right)
,
\left(\begin{array}{ccccc}
0&0&r_{12}&0\\
0&0&0&0\\
r_{31}&r_{32}&r_{33}&0\\
r_{41}&r_{42}&r_{43}&0\\
\end{array}
\right)
,
\left(\begin{array}{ccccc}
0&0&0&0\\
0&0&0&0\\
r_{31}&r_{32}&r_{33}&r_{34}\\
r_{41}&r_{42}&r_{43}&r_{44}\\
\end{array}
\right)
;
$

\item[$L_{20}$ : ]
    $\left(\begin{array}{ccccc}
0&0&0&0\\
r_{21}&0&0&0\\
r_{31}&r_{32}&0&\frac{1}{2}(r_{21}+\alpha r_{21})\\
r_{41}&r_{42}&0&0 
\end{array}
\right),\quad
\left(\begin{array}{ccccc}
r_{11}&r_{12}&0&0\\
0&0&0&0\\
r_{31}&r_{32}&r_{33}&0\\
r_{41}&r_{42}&r_{43}&0
\end{array}
\right)
,\quad
\left(\begin{array}{ccccc}
r_{11}&r_{12}&r_{13}&0\\
0&0&0&0\\
r_{31}&r_{32}&r_{33}&0\\
r_{41}&r_{42}&r_{43}&0\\
\end{array}
\right)
,
$

 $\left(\begin{array}{ccccc}
 \frac{r_{22}r_{44}}{r_{22}+r_{44}} &0&0&0\\
0&r_{22}&0&0\\
r_{31}&r_{32}&\frac{1}{2}r_{22}&0\\
r_{41}&r_{42}&0&r_{44}\\
\end{array}
\right)
,\quad
\left(\begin{array}{ccccc}
0&0&r_{12}&0\\
0&0&0&0\\
r_{31}&r_{32}&r_{33}&0\\
r_{41}&r_{42}&r_{43}&0\\
\end{array}
\right)
,
\quad
\left(\begin{array}{ccccc}
0&r_{12}&0&0\\
0&r_{22}&0&0\\
r_{31}&r_{32}&\frac{1}{2}r_{22}&0\\
r_{41}&r_{42}&\frac{r_{22}}{1-\alpha}&0\\
\end{array}
\right)
;
$

\item[$L_{21}$ : ]
    $\left(\begin{array}{ccccc}
r_{11}&r_{12}&r_{13}&0\\
0&0&0&0\\
0&0&0&0\\
r_{41}&r_{42}&r_{43}&0 
\end{array}
\right),\quad
\left(\begin{array}{ccccc}
0&0&0&0\\
0&0&0&0\\
0&0&0&0\\
r_{41}&r_{42}&r_{43}&r_{44}
\end{array}
\right)
,\quad
\left(\begin{array}{ccccc}
0&0&0&0\\
0&0&0&0\\
0&0&2r_{44}&0\\
r_{41}&r_{42}&r_{43}&r_{44}
\end{array}
\right)
,
$

 $\left(\begin{array}{ccccc}
 \frac{1}{2}r_{22}&0&0&0\\
0&r_{22}&0&0\\
0&0&\frac{2}{3}r_{22}&0\\
r_{41}&r_{42}&0&\frac{1}{3}r_{22}\\
\end{array}
\right)
,\quad
\left(\begin{array}{ccccc}
0&r_{22}&0&0\\
0&0&0&0\\
r_{41}&r_{42}&0&\frac{1}{3}r_{22}\\
\end{array}
\right)
,
\quad
\left(\begin{array}{ccccc}
0&r_{12}&r_{13}&0\\
0&0&0&0\\
0&0&0&0\\
r_{41}&r_{42}&r_{43}&0\\
\end{array}
\right).
$
\end{itemize}
\end{prop}

\subsection{Nenjenhuis operator}\,

Let $(L, \llbracket-,-\rrbracket)$   be an $n$-dimensional Leibniz algebra with basis $\{e_i\}\, (1 \leq i \leq n)$ and let $\kappa$ be a Nenjenhuis operator on $L$
For any $i,j,k\in\mathbb{N}$, $1 \leq i,j,k \leq n$, let us put
$$\llbracket e_i,e_j\rrbracket=\sum_{k=1}^nc_{ij}^ke_k\quad \text{and}\quad \kappa(e_i)=\sum_{j=1}^n\kappa_{ji}e_j.$$
Then, in term of basis elements, equation (\ref{eq4}) is equivalent to
\begin{align*}
   \sum^{n}_{p=1}\sum^{n}_{q=1}\kappa_{pi}\kappa_{qj} C_{kq}^{t}= 
   \sum^{n}_{p=1} \sum^{n}_{q=1}\Bigg(\kappa_{pi}C_{pj}^{q}\kappa_{tq}+\kappa_{pj}C_{ip}^{q}\kappa_{tq}-c_{ij}^p\kappa_{qp}\kappa_{tq}\Bigg), 
\end{align*} for $i,j,t=1,2,\dots,n$.

\begin{prop}
The description of the Nenjenhuis operator of every 4-dimensional  Leibniz algebra is given below.
\begin{itemize}
    \item[$L_1$ : ]
    $\left(\begin{array}{ccccc}
\kappa_{11}&0&0&0\\
\kappa_{21}&\kappa_{11}&0&0\\
\kappa_{31}&\kappa_{32}&\kappa_{11}&0\\
\kappa_{41}&\kappa_{42}&\kappa_{43}&\kappa_{11}\\
\end{array}
\right),\,
\left(\begin{array}{ccccc}
\kappa_{11}&0&0&0\\
\kappa_{21}&\kappa_{11}&0&0\\
\kappa_{31}&\kappa_{32}&\kappa_{11}&0\\
\kappa_{41}&\kappa_{42}&0&\kappa_{11}\\
\end{array}
\right)
,
\left(\begin{array}{ccccc}
\kappa_{11}&0&0&0\\
\kappa_{21}&\kappa_{11}&0&0\\
0&0&\kappa_{11}&\kappa_{24}\\
\kappa_{41}&-\frac{\kappa_{11}-\kappa_{22}^2}{\kappa_{24}}&0&2\kappa_{11}-\kappa_{22}\\
\end{array}
\right)
;
$

\item[$L_2$ : ]
    $\left(\begin{array}{ccccc}
\kappa_{11}&0&0&0\\
\kappa_{21}&\kappa_{22}&0&0\\
\kappa_{31}&\kappa_{32}&\kappa_{11}&0\\
\kappa_{41}&\kappa_{42}&0&\kappa_{11}\\
\end{array}
\right),\,
\left(\begin{array}{ccccc}
\kappa_{11}&0&0&0\\
0&\kappa_{11}&0&0\\
\kappa_{31}&\kappa_{32}&\kappa_{11}&0\\
\kappa_{41}&\kappa_{42}&\kappa_{43}&\kappa_{11}\\
\end{array}
\right)
,\,
\left(\begin{array}{ccccc}
\kappa_{11}&0&0&0\\
\kappa_{21}&\kappa_{11}&0&0\\
\kappa_{31}&0&\kappa_{11}&0\\
\kappa_{41}&\kappa_{42}&0&\kappa_{11}\\
\end{array}
\right),\\
\left(\begin{array}{ccccc}
\kappa_{11}&0&0&0\\
\kappa_{21}&\kappa_{11}&0&0\\
\kappa_{31}&0&\kappa_{11}&0\\
\kappa_{41}&\kappa_{42}&\kappa_{43}&\kappa_{11}\\
\end{array}
\right),\quad
\left(\begin{array}{ccccc}
\kappa_{11}&0&0&0\\
\kappa_{33}-\kappa_{11}&\kappa_{33}&0&0\\
\kappa_{31}&\kappa_{32}&\kappa_{33}&0\\
\kappa_{41}&\kappa_{42}&\kappa_{43}&\kappa_{11}\\
\end{array}
\right)
;
$
\item[$L_3$ : ]
    $\left(\begin{array}{ccccc}
\kappa_{11}&0&0&0\\
\kappa_{21}&\kappa_{22}&0&0\\
\kappa_{31}&\kappa_{32}&\kappa_{11}&0\\
\kappa_{41}&\kappa_{42}&0&\kappa_{11}\\
\end{array}
\right),\quad
\left(\begin{array}{ccccc}
\kappa_{11}&0&0&0\\
0&\kappa_{11}&0&0\\
\kappa_{31}&\kappa_{32}&\kappa_{11}&0\\
\kappa_{41}&\kappa_{42}&\kappa_{43}&\kappa_{11}\\
\end{array}
\right)
,\quad
\left(\begin{array}{ccccc}
\kappa_{11}&0&0&0\\
\kappa_{33}-\kappa_{21}&\kappa_{33}&0&0\\
\kappa_{31}&\kappa_{32}&\kappa_{33}&0\\
\kappa_{41}&\kappa_{42}&\kappa_{43}&\kappa_{11}\\
\end{array}
\right)
;
$

\item[$L_4$ : ]
    $\left(\begin{array}{ccccc}
\kappa_{44}&0&0&0\\
\kappa_{21}&\kappa_{22}&0&0\\
\kappa_{31}&\kappa_{32}&\kappa_{44}&0\\
0&0&0&\kappa_{44}\\
\end{array}
\right),\quad
\left(\begin{array}{ccccc}
\kappa_{44}&0&0&0\\
0&\kappa_{44}&0&0\\
\kappa_{31}&\kappa_{32}&\kappa_{44}&0\\
\kappa_{41}&\kappa_{42}&\kappa_{43}&\kappa_{44}\\
\end{array}
\right)
,\quad
\left(\begin{array}{ccccc}
\kappa_{44}&0&0&0\\
\kappa_{21}&\kappa_{22}&0&0\\
\kappa_{31}&\kappa_{32}&\kappa_{44}&0\\
0&0&0&\kappa_{44}\\
\end{array}
\right),
$
\\

$
\left(\begin{array}{ccccc}
\kappa_{44}&0&0&0\\
0&\kappa_{22}&0&0\\
\kappa_{31}&\kappa_{32}&\kappa_{44}&0\\
0&0&0&\kappa_{44}\\
\end{array}
\right)
,\quad
\left(\begin{array}{ccccc}
\kappa_{44}&0&0&0\\
0&\kappa_{22}&0&0\\
\kappa_{31}&\kappa_{32}&\kappa_{44}&0\\
0&0&0&\kappa_{44}\\
\end{array}
\right)
;
$

\item[$L_5$ : ]
    $\left(\begin{array}{ccccc}
\kappa_{44}&0&0&0\\
0&\kappa_{44}&0&0\\
\kappa_{31}&\kappa_{32}&\kappa_{44}&0\\
\kappa_{41}&\kappa_{42}&\kappa_{43}&\kappa_{44}\\
\end{array}
\right),
\left(\begin{array}{ccccc}
\kappa_{44}&0&0&0\\
0&\kappa_{44}&0&0\\
\kappa_{31}&\kappa_{32}&\kappa_{44}&0\\
\kappa_{41}&\kappa_{42}&0&\kappa_{44}\\
\end{array}
\right)
,
\left(\begin{array}{ccccc}
\kappa_{44}&0&0&0\\
\kappa_{21}&\kappa_{22}&0&0\\
\kappa_{31}&\kappa_{32}&\kappa_{44}&0\\
0&0&0&\kappa_{44}\\
\end{array}
\right),
\left(\begin{array}{ccccc}
\kappa_{44}&0&0&0\\
0&\kappa_{22}&0&0\\
\kappa_{31}&\kappa_{32}&\kappa_{44}&0\\
0&0&0&\kappa_{44}\\
\end{array}
\right)
;
$

\item[$L_6$ : ]
    $\left(\begin{array}{ccccc}
\kappa_{44}&0&0&0\\
0&\kappa_{22}&0&0\\
\kappa_{31}&\kappa_{32}&\kappa_{44}&0\\
\kappa_{41}&\kappa_{42}&\kappa_{43}&\kappa_{44}\\
\end{array}
\right)
;
$

\item[$L_7$ : ]
    $
\left(\begin{array}{ccccc}
\kappa_{11}&0&0&0\\
\kappa_{21}&\kappa_{22}&0&0\\
\kappa_{31}&\kappa_{32}&\kappa_{33}&0\\
\kappa_{41}&\kappa_{42}&0&\kappa_{11}\\
\end{array}
\right)
,\quad
\left(\begin{array}{ccccc}
\kappa_{11}&0&0&0\\
\kappa_{21}&\kappa_{22}&0&0\\
\kappa_{31}&\kappa_{32}&\kappa_{11}&0\\
\kappa_{41}&\kappa_{42}&0&\kappa_{11}\\
\end{array}
\right),
,\quad
\left(\begin{array}{ccccc}
\kappa_{11}&0&0&0\\
\kappa_{21}&\kappa_{22}&0&0\\
\kappa_{31}&\kappa_{32}&\kappa_{11}&0\\
\kappa_{31}&\kappa_{32}&\kappa_{43}&\kappa_{11}\\
\end{array}
\right)
;
$

\item[$L_8$ : ]
    $
\left(\begin{array}{ccccc}
\kappa_{44}&0&0&0\\
\kappa_{21}&\kappa_{33}&0&0\\
\kappa_{31}&\kappa_{32}&\kappa_{33}&0\\
\kappa_{41}&\kappa_{42}&\kappa_{43}&\kappa_{11}\\
\end{array}
\right)
,
\left(\begin{array}{ccccc}
\kappa_{44}&0&0&0\\
\kappa_{21}&\kappa_{44}&0&0\\
\kappa_{31}&\kappa_{32}&\kappa_{44}&0\\
\kappa_{41}&\kappa_{42}&\kappa_{43}&\kappa_{44}\\
\end{array}
\right),
\left(\begin{array}{ccccc}
\kappa_{44}&0&0&0\\
\kappa_{21}&\kappa_{22}&0&0\\
\kappa_{31}&\kappa_{32}&\kappa_{44}&0\\
\kappa_{31}&\kappa_{32}&0&\kappa_{44}\\
\end{array}
\right)
,
\left(\begin{array}{ccccc}
\kappa_{44}&0&0&0\\
\kappa_{21}&\kappa_{33}&0&0\\
\kappa_{31}&\kappa_{32}&\kappa_{33}&0\\
\kappa_{31}&\kappa_{32}&0&\kappa_{44}\\
\end{array}
\right)
;
$

\item[$L_9$ : ]
    $
\left(\begin{array}{ccccc}
\kappa_{44}&\kappa_{21}&0&0\\
0&\kappa_{44}-\kappa_{11}&0&0\\
\kappa_{31}&\kappa_{32}&\kappa_{44}&0\\
\kappa_{41}&\kappa_{42}&\kappa_{43}&\kappa_{44}\\
\end{array}
\right)
,\,
\left(\begin{array}{ccccc}
\kappa_{44}&0&0&0\\
\kappa_{21}&\kappa_{44}&0&0\\
\kappa_{31}&\kappa_{32}&\kappa_{44}&0\\
\kappa_{41}&\kappa_{42}&\kappa_{43}&\kappa_{44}\\
\end{array}
\right),
\,
\left(\begin{array}{ccccc}
\kappa_{44}&0&0&0\\
\kappa_{44}-\kappa_{11}&\kappa_{22}&0&0\\
\kappa_{31}&\kappa_{32}&\kappa_{44}&0\\
\kappa_{31}&\kappa_{32}&\kappa_{43}&\kappa_{44}\\
\end{array}
\right)
;
$

\item[$L_{10}$ : ]
    $
\left(\begin{array}{ccccc}
\kappa_{44}&\kappa_{21}&0&0\\
0&\kappa_{44}&0&0\\
\kappa_{31}&\kappa_{32}&\kappa_{44}&0\\
\kappa_{41}&\kappa_{42}&\kappa_{43}&\kappa_{44}\\
\end{array}
\right)
,\,
\left(\begin{array}{ccccc}
\kappa_{44}&\kappa_{12}&0&0\\
0&\kappa_{44}&0&0\\
0&-\frac{1}{2}\kappa_{12}&\kappa_{44}&0\\
\kappa_{41}&\kappa_{42}&\kappa_{43}&\kappa_{44}\\
\end{array}
\right),
\,
\left(\begin{array}{ccccc}
\kappa_{11}&\kappa_{12}&0&0\\
0&\kappa_{22}&0&0\\
0&-\frac{1}{2}\kappa_{12}&\frac{\kappa_{44}+\kappa_{11}}{2}&0\\
\kappa_{31}&\kappa_{32}&0&\kappa_{44}\\
\end{array}
\right)
;
$

\item[$L_{11}$ : ]
    $
\left(\begin{array}{ccccc}
\kappa_{11}&0&0&0\\
0&\kappa_{33}&0&0\\
\kappa_{33}-\kappa_{11}&0&\kappa_{33}&0\\
\kappa_{41}&\kappa_{42}&0&\kappa_{44}\\
\end{array}
\right)
,\,
\left(\begin{array}{ccccc}
\kappa_{33}&0&0&0\\
0&\kappa_{33}&0&0\\
\kappa_{31}&\kappa_{32}&\kappa_{33}&0\\
\kappa_{31}&\kappa_{32}&\kappa_{43}&\kappa_{33}\\
\end{array}
\right),
\left(\begin{array}{ccccc}
\kappa_{33}&\kappa_{12}&0&0\\
0&\kappa_{11}+\kappa_{33}&0&0\\
0&2\kappa_{12}&\kappa_{33}&0\\
\kappa_{41}&\kappa_{42}&0&\kappa_{33}\\
\end{array}
\right)
;
$

\item[$L_{12}$ : ]
    $
\left(\begin{array}{ccccc}
\kappa_{11}&\kappa_{12}&0&0\\
0&\kappa_{22}&0&0\\
\kappa_{31}&\kappa_{32}&\kappa_{22}&0\\
\kappa_{41}&\kappa_{42}&0&\kappa_{22}\\
\end{array}
\right)
,
\left(\begin{array}{ccccc}
\kappa_{11}&0&0&0\\
0&\kappa_{22}&0&0\\
\kappa_{31}&\kappa_{32}&\kappa_{22}&0\\
\kappa_{41}&\kappa_{42}&0&\kappa_{22}\\
\end{array}
\right),
\left(\begin{array}{ccccc}
\kappa_{44}&0&0&0\\
0&\kappa_{22}&0&0\\
\kappa_{31}&\kappa_{32}&\kappa_{22}&0\\
\kappa_{31}&\kappa_{32}&\kappa_{43}&\kappa_{22}\\
\end{array}
\right),
\left(\begin{array}{ccccc}
\kappa_{44}&0&0&0\\
0&\kappa_{22}&0&0\\
\kappa_{31}&\kappa_{32}&\kappa_{22}&0\\
\kappa_{31}&\kappa_{32}&0&\kappa_{44}\\
\end{array}
\right)
;
$

\item[$L_{13}$ : ]
    $
\left(\begin{array}{ccccc}
\kappa_{33}&0&0&0\\
0&\kappa_{33}&0&0\\
\kappa_{31}&\kappa_{32}&\kappa_{33}&0\\
\kappa_{41}&\kappa_{42}&\kappa_{43}&\kappa_{33}\\
\end{array}
\right)
,
\left(\begin{array}{ccccc}
\kappa_{33}&0&0&0\\
0&\kappa_{22}&0&0\\
\kappa_{31}&\kappa_{32}&\kappa_{33}&0\\
\kappa_{41}&\kappa_{42}&0&\kappa_{22}\\
\end{array}
\right),
\left(\begin{array}{ccccc}
\kappa_{33}&0&0&0\\
0&\kappa_{33}&0&0\\
\kappa_{31}&\kappa_{32}&\kappa_{33}&\kappa_{34}\\
\kappa_{41}&\kappa_{42}&0&\kappa_{33}\\\
\end{array}
\right),
\left(\begin{array}{ccccc}
\kappa_{33}&0&0&0\\
0&\kappa_{33}&0&0\\
\kappa_{31}&\kappa_{32}&\kappa_{33}&0\\
\kappa_{41}&\kappa_{42}&0&N_{33}\\
\end{array}
\right)
;
$

\item[$L_{14}$ : ]
    $
\left(\begin{array}{ccccc}
\kappa_{44}&0&0&0\\
\kappa_{21}&\kappa_{22}&\kappa_{23}&0\\
-i\kappa_{21}&i(\kappa_{44}-\kappa_{22})&\kappa_{44}&0\\
\kappa_{41}&\kappa_{42}&\kappa_{43}&\kappa_{44}\\
\end{array}
\right)
,
\left(\begin{array}{ccccc}
\kappa_{11}&\kappa_{12}&0&0\\
0&\kappa_{44}&0&0\\
i(\kappa_{44}-\kappa_{22})&-\kappa_{12}&\kappa_{44}&0\\
\kappa_{41}&\kappa_{42}&\kappa_{43}&\kappa_{44}\\
\end{array}
\right),
\left(\begin{array}{ccccc}
\kappa_{11}&0&0&0\\
i(\kappa_{44}-\kappa_{22})&\kappa_{44}&0&0\\
0&0&\kappa_{44}&0\\
\kappa_{41}&\kappa_{42}&\kappa_{43}&\kappa_{44}\\
\end{array}
\right)
;
$

\item[$L_{15}$ : ]
    $
\left(\begin{array}{ccccc}
\kappa_{44}&0&0&0\\
0&\kappa_{44}&0&0\\
\kappa_{31}&\kappa_{32}&\kappa_{33}&0\\
\kappa_{41}&\kappa_{42}&\kappa_{43}&\kappa_{44}\\
\end{array}
\right)
,\,
\left(\begin{array}{ccccc}
\kappa_{11}&\kappa_{12}&\kappa_{13}&0\\
0&\kappa_{44}&0&0\\
0&0&\kappa_{44}&0\\
\kappa_{41}&\kappa_{42}&\kappa_{43}&\kappa_{44}\\
\end{array}
\right),
\,
\left(\begin{array}{ccccc}
\kappa_{44}&0&0&0\\
0&\kappa_{44}&0&0\\
0&\kappa_{32}&\kappa_{33}&0\\
\kappa_{41}&\kappa_{42}&\kappa_{43}&\kappa_{44}\\
\end{array}
\right)
,
\left(\begin{array}{ccccc}
\kappa_{11}&0&0&0\\
0&\kappa_{44}&0&0\\
0&0&\kappa_{44}&0\\
\kappa_{41}&\kappa_{42}&\kappa_{43}&\kappa_{44}\\
\end{array}
\right)
$

\item[$L_{16}$ : ]
    $
\left(\begin{array}{ccccc}
\kappa_{44}&0&0&0\\
\kappa_{21}&\kappa_{22}&\kappa_{23}&0\\
0&0&\kappa_{44}&0\\
\kappa_{41}&\kappa_{42}&\kappa_{43}&\kappa_{44}\\
\end{array}
\right)
,\,
\left(\begin{array}{ccccc}
\kappa_{44}&0&0&0\\
\kappa_{21}&\kappa_{22}&0&0\\
0&0&\kappa_{44}&0\\
\kappa_{41}&\kappa_{42}&\kappa_{43}&\kappa_{44}\\
\end{array}
\right),
$

\item[$L_{17}$ : ]
    $
\left(\begin{array}{ccccc}
\kappa_{11}&0&0&0\\
0&\kappa_{22}&0&0\\
\kappa_{31}&\kappa_{32}&\kappa_{22}&0\\
\kappa_{41}&\kappa_{42}&0&\kappa_{22}\\
\end{array}
\right)
,\,
\left(\begin{array}{ccccc}
\kappa_{11}&0&0&0\\
0&\kappa_{22}&0&0\\
\kappa_{31}&\kappa_{32}&\kappa_{22}&0\\
\kappa_{41}&\kappa_{42}&0&\kappa_{11}\\
\end{array}
\right),
\,
\left(\begin{array}{ccccc}
\kappa_{11}&0&0&0\\
0&\kappa_{22}&0&0\\
\kappa_{31}&\kappa_{32}&\kappa_{11}&0\\
\kappa_{41}&\kappa_{42}&0&\kappa_{11}\\
\end{array}
\right)
;
$

$
\left(\begin{array}{ccccc}
\kappa_{11}&0&0&0\\
\kappa_{21}&\kappa_{11}&0&0\\
\kappa_{31}&\kappa_{32}&\kappa_{11}&0\\
\kappa_{41}&\kappa_{42}&0&\kappa_{11}\\
\end{array}
\right),\quad
\left(\begin{array}{ccccc}
\kappa_{11}&0&0&0\\
0&\kappa_{22}&0&0\\
\kappa_{31}&\kappa_{32}&\kappa_{22}&0\\
\kappa_{41}&\kappa_{42}&\kappa_{43}&\kappa_{11}\\
\end{array}
\right);
$

\item[$L_{18}$ : ]
    $
\left(\begin{array}{ccccc}
\kappa_{33}&\frac{\kappa_{34}}{2}&0&0\\
2\kappa_{43}&\kappa_{22}&0&0\\
\kappa_{31}&\kappa_{32}&\kappa_{33}&\kappa_{34}\\
\kappa_{41}&\kappa_{42}&0&\kappa_{22}\\
\end{array}
\right)
,\,
\left(\begin{array}{ccccc}
\kappa_{11}&\kappa_{12}&0&0\\
0&\kappa_{22}&0&0\\
\kappa_{31}&\kappa_{32}&\kappa_{22}&\kappa_{34}\\
\kappa_{41}&\kappa_{42}&0&\kappa_{22}\\
\end{array}
\right),
\,
\left(\begin{array}{ccccc}
\kappa_{33}&0&0&0\\
0&\kappa_{22}&0&0\\
\kappa_{31}&\kappa_{32}&\kappa_{33}&0\\
\kappa_{41}&\kappa_{42}&\kappa_{43}&\kappa_{22}\\
\end{array}
\right)
;
$

$
\left(\begin{array}{ccccc}
\kappa_{11}&\kappa_{12}&0&0\\
0&\kappa_{22}&0&0\\
\kappa_{31}&\kappa_{32}&\kappa_{33}&0\\
\kappa_{41}&\kappa_{42}&0&\kappa_{22}\\
\end{array}
\right),\quad
\left(\begin{array}{ccccc}
\kappa_{33}&0&0&0\\
0&\kappa_{22}&0&0\\
\kappa_{31}&\kappa_{32}&\kappa_{33}&0\\
\kappa_{41}&\kappa_{42}&\kappa_{43}&\kappa_{22}\\
\end{array}
\right);
$

\item[$L_{19}$ : ]
    $
\left(\begin{array}{ccccc}
\kappa_{11}&\kappa_{12}&0&0\\
0&\kappa_{22}&0&0\\
\kappa_{31}&\kappa_{32}&\kappa_{33}&0\\
\kappa_{41}&\kappa_{42}&0&\kappa_{22}\\
\end{array}
\right)
,
\left(\begin{array}{ccccc}
\kappa_{11}&\kappa_{12}&0&0\\
0&\kappa_{22}&0&0\\
\kappa_{31}&\kappa_{32}&\kappa_{22}&0\\
\kappa_{41}&v_{42}&0&\kappa_{44}\\
\end{array}
\right),
\,
\left(\begin{array}{ccccc}
\kappa_{44}&0&0&0\\
0&v_{22}&0&0\\
\kappa_{31}&\kappa_{32}&\kappa_{33}&\kappa_{34}\\
\kappa_{41}&\kappa_{42}&0&\kappa_{22}\\
\end{array}
\right)
,
\left(\begin{array}{ccccc}
\kappa_{22}&\kappa_{12}&\kappa_{13}&0\\
0&\kappa_{22}&0&0\\
0&\kappa_{32}&\kappa_{22}&0\\
0&\kappa_{42}&0&\kappa_{22}\\
\end{array}
\right);
$

\item[$L_{20}$ : ]
    $
\left(\begin{array}{ccccc}
\kappa_{11}&\kappa_{12}&\kappa_{13}&0\\
0&\kappa_{44}&0&0\\
\kappa_{31}&\kappa_{32}&\kappa_{33}&0\\
\kappa_{41}&\kappa_{42}&\kappa_{43}&\kappa_{44}\\
\end{array}
\right)
,
\left(\begin{array}{ccccc}
\kappa_{44}&0&0&0\\
0&\kappa_{22}&\kappa_{23}&0\\
\kappa_{31}&\kappa_{32}&\kappa_{33}&0\\
\kappa_{41}&\kappa_{42}&\kappa_{43}&\kappa_{44}\\
\end{array}
\right),
\left(\begin{array}{ccccc}
\kappa_{44}&0&0&0\\
\kappa_{21}&\kappa_{44}&0&0\\
\kappa_{31}&\kappa_{32}&\kappa_{33}&0\\
\kappa_{41}&\kappa_{42}&\kappa_{43}&\kappa_{44}
\end{array}
\right),
\left(\begin{array}{ccccc}
\kappa_{11}&0&\kappa_{13}&0\\
0&\kappa_{44}&0&0\\
0&0&\kappa_{33}&0\\
\kappa_{41}&\kappa_{42}&\kappa_{43}&\kappa_{44}\\
\end{array}
\right);
$

\item[$L_{21}$ : ]
    $
\left(\begin{array}{ccccc}
\kappa_{44}&0&0&0\\
\kappa_{21}&\kappa_{22}&\kappa_{23}&0\\
0&0&\kappa_{33}&0\\
\kappa_{41}&\kappa_{42}&\kappa_{43}&\kappa_{44}\\
\end{array}
\right)
,
\left(\begin{array}{ccccc}
\kappa_{11}&\kappa_{12}&\kappa_{13}&0\\
0&\kappa_{44}&0&0\\
0&0&\kappa_{44}&0\\
\kappa_{41}&\kappa_{42}&\kappa_{43}&\kappa_{44}\\
\end{array}
\right),
\left(\begin{array}{ccccc}
\kappa_{44}&0&0&0\\
\kappa_{21}&\kappa_{22}&0&0\\
0&0&\kappa_{44}&0\\
\kappa_{41}&\kappa_{42}&\kappa_{43}&\kappa_{44}\\
\end{array}
\right)
,
\left(\begin{array}{ccccc}
\kappa_{11}&0&\kappa_{13}&0\\
0&\kappa_{44}&0&0\\
0&0&\kappa_{33}&0\\
\kappa_{41}&\kappa_{42}&\kappa_{43}&\kappa_{44}\\
\end{array}
\right).
$
\end{itemize}
\end{prop}

\subsection{Reynolds operator}\,

Let $(L, \llbracket-,-\rrbracket)$   be an $n$-dimensional Leibniz algebra with basis $\{e_i\}\, (1 \leq i \leq n)$ and let $\alpha$ be a Reynolds operator on $L$.
For any $i,j,k\in\mathbb{N}$, $1 \leq i,j,k \leq n$, let us put
$$\llbracket e_i,e_j\rrbracket=\sum_{k=1}^nc_{ij}^ke_k\quad \text{and}\quad \alpha(e_i)=\sum_{j=1}^n\alpha_{ji}e_j.$$
Then, in term of basis elements, equation (\ref{eq5}) is equivalent to
\begin{align*}
   \sum^{n}_{p=1}\sum^{n}_{q=1}\alpha_{pi}\alpha_{qj} C_{kq}^{t}= 
   \sum^{n}_{p=1} \sum^{n}_{q=1}\Bigg(\alpha_{pi}C_{pj}^{q}\alpha_{tq}+\alpha_{pj}C_{ip}^{q}\alpha_{tq}\Bigg)-\sum_{q=1}^n\sum_{k=1}^n\sum_{p=1}^n\alpha_{ki}\alpha_{pj}C_{kp}^q\alpha_{tq}, 
\end{align*} for $i,j,t=1,2,\dots,n$.

\begin{prop}
The description of the Reynolds operator of every 4-dimensional  Leibniz algebra is given below.
\begin{itemize}
    \item[$L_1$ : ]
    $\left(\begin{array}{ccccc}
\alpha_{11}&0&0&0\\
0&\alpha_{11}&0&0\\
0&0&\alpha_{11}&0\\
\alpha_{41}&0&0&\alpha_{11}\\
\end{array}
\right),\,
\left(\begin{array}{ccccc}
0&0&0&0\\
0&0&0&0\\
\alpha_{31}&\alpha_{32}&\alpha_{11}&0\\
\alpha_{41}&\alpha_{42}&\alpha_{43}&0\\
\end{array}
\right)
,\,
\left(\begin{array}{ccccc}
0&0&0&0\\
\alpha_{21}&\alpha_{22}&0&\alpha_{24}\\
0&0&0&0\\
\alpha_{41}&\alpha_{42}&0&\alpha_{44}\\
\end{array}
\right)
,\,
\left(\begin{array}{ccccc}
0&0&0&0\\
0&0&0&0\\
\alpha_{21}&0&\alpha_{33}&0\\
\alpha_{41}&\alpha_{42}&\alpha_{43}&0\\
\end{array}
\right)
;
$

\item[$L_2$ : ]
     $\left(\begin{array}{ccccc}
\alpha_{11}&0&0&0\\
0&\alpha_{11}&0&0\\
0&0&\alpha_{11}&0\\
\alpha_{41}&\alpha_{42}&0&\alpha_{11}\\
\end{array}
\right),\,
\left(\begin{array}{ccccc}
0&0&0&0\\
0&0&0&0\\
\alpha_{31}&\alpha_{32}&\alpha_{33}&0\\
\alpha_{41}&\alpha_{42}&\alpha_{43}&0\\
\end{array}
\right)
,\,
\left(\begin{array}{ccccc}
0&0&0&0\\
\alpha_{21}&\alpha_{22}&0&0\\
\alpha_{31}&\alpha_{32}&0&0\\
\alpha_{41}&\alpha_{42}&0&0\\
\end{array}
\right)
,\,
\left(\begin{array}{ccccc}
0&0&0&0\\
0&0&0&0\\
\alpha_{31}&\alpha_{32}&0&0\\
\alpha_{41}&\alpha_{42}&\alpha_{43}&0\\
\end{array}
\right)
;
$
\item[$L_3$ : ]
      $\left(\begin{array}{ccccc}
\alpha_{11}&0&0&0\\
0&\alpha_{11}&0&0\\
0&0&\alpha_{11}&0\\
\alpha_{41}&\alpha_{42}&0&\alpha_{11}\\
\end{array}
\right),\,
\left(\begin{array}{ccccc}
0&0&0&0\\
0&0&0&0\\
\alpha_{31}&\alpha_{32}&\alpha_{33}&0\\
\alpha_{41}&\alpha_{42}&\alpha_{43}&0\\
\end{array}
\right)
,\,
\left(\begin{array}{ccccc}
0&0&0&0\\
\alpha_{21}&\alpha_{22}&0&0\\
0&0&0&0\\
\alpha_{41}&\alpha_{42}&0&\alpha_{44}\\
\end{array}
\right)
,\,
\left(\begin{array}{ccccc}
0&0&0&0\\
\alpha_{31}&\alpha_{32}&0&\alpha_{24}\\
0&0&0&0\\
\alpha_{41}&\alpha_{42}&0&\alpha_{44}\\
\end{array}
\right)
;
$

\item[$L_4$ : ]
    $\left(\begin{array}{ccccc}
\alpha_{11}&0&0&0\\
0&\alpha_{11}&0&0\\
0&0&\alpha_{11}&0\\
\alpha_{41}&\alpha_{42}&0&\alpha_{11}\\
\end{array}
\right),\,
\left(\begin{array}{ccccc}
0&0&0&0\\
0&0&0&0\\
\alpha_{31}&\alpha_{32}&\alpha_{33}&0\\
\alpha_{41}&\alpha_{42}&\alpha_{43}&0\\
\end{array}
\right)
,\,
\left(\begin{array}{ccccc}
0&0&0&0\\
0&0&0&0\\
0&0&0&0\\
\alpha_{41}&\alpha_{42}&\alpha_{43}&\alpha_{44}\\
\end{array}
\right)
;
$

\item[$L_5$ : ]
    $\left(\begin{array}{ccccc}
\alpha_{11}&0&0&0\\
\alpha_{21}&\alpha_{11}&0&0\\
\alpha_{31}&0&\alpha_{11}&0\\
\alpha_{41}&\alpha_{42}&\alpha_{31}&\alpha_{11}\\
\end{array}
\right),\,
\left(\begin{array}{ccccc}
0&0&0&0\\
0&0&0&0\\
0&0&0&0\\
\alpha_{41}&\alpha_{42}&\alpha_{43}&\alpha_{44}\\
\end{array}
\right)
,\,
\left(\begin{array}{ccccc}
0&0&0&0\\
\alpha_{21}&\alpha_{22}&\alpha_{23}&0\\
\alpha_{31}&\alpha_{32}&\alpha_{23}&0\\
\alpha_{41}&\alpha_{42}&0&\alpha_{44}\\
\end{array}
\right)
,\,
\left(\begin{array}{ccccc}
0&0&0&0\\
\alpha_{21}&\alpha_{22}&0&0\\
\alpha_{31}&\alpha_{32}&\alpha_{33}&0\\
\alpha_{41}&\alpha_{42}&\alpha_{34}&0\\
\end{array}
\right)
;
$

\item[$L_6$ : ]
    $\left(\begin{array}{ccccc}
0&0&0&0\\
0&0&0&0\\
\alpha_{31}&\alpha_{32}&\alpha_{33}&0\\
\alpha_{41}&\alpha_{42}&\alpha_{43}&0\\
\end{array}
\right),\,
\left(\begin{array}{ccccc}
0&0&0&0\\
0&0&0&0\\
\alpha_{31}&\alpha_{32}&0&0\\
\alpha_{41}&\alpha_{42}&\alpha_{43}&0\\
\end{array}
\right)
,\,
\left(\begin{array}{ccccc}
\alpha_{11}&0&0&0\\
0&0&0&0\\
0&0&\alpha_{11}&0\\
\alpha_{41}&\alpha_{42}&0&\alpha_{11}\\
\end{array}
\right)
,\,
\left(\begin{array}{ccccc}
\alpha_{11}&0&0&0\\
0&\alpha_{11}&0&0\\
0&0&\alpha_{11}&0\\
\alpha_{41}&\alpha_{42}&0&\alpha_{11}\\
\end{array}
\right)
,$
\\
$\left(\begin{array}{ccccc}
0&0&0&0\\
0&0&0&0\\
0&0&0&0\\
\alpha_{41}&\alpha_{42}&\alpha_{43}&\alpha_{44}\\
\end{array}
\right)
,\,
\left(\begin{array}{ccccc}
0&0&0&0\\
0&\alpha_{22}&0&0\\
0&0&0&0\\
\alpha_{41}&\alpha_{42}&\alpha_{43}&\alpha_{22}\\
\end{array}
\right)$

\item[$L_7$ : ]
    $\left(\begin{array}{ccccc}
0&0&0&0\\
0&\alpha_{22}&0&0\\
0&0&0&0\\
\alpha_{41}&\alpha_{42}&0&\alpha_{44}\\
\end{array}
\right),\,
\left(\begin{array}{ccccc}
0&0&0&0\\
\alpha_{21}&\alpha_{22}&0&0\\
0&0&0&0\\
\alpha_{41}&\alpha_{42}&0&\alpha_{44}\\
\end{array}
\right)
,\,
\left(\begin{array}{ccccc}
0&0&0&0\\
0&\alpha_{22}&0&\alpha_{24}\\
0&0&0&0\\
\alpha_{41}&\alpha_{42}&0&\alpha_{44}\\
\end{array}
\right)
,\,
\left(\begin{array}{ccccc}
\alpha_{11}&0&0&0\\
0&\alpha_{11}&0&0\\
0&0&\alpha_{11}&0\\
\alpha_{41}&\alpha_{42}&0&\alpha_{11}\\
\end{array}
\right)
,$
\\
$\left(\begin{array}{ccccc}
0&0&0&0\\
0&0&0&0\\
\alpha_{31}&\alpha_{32}&\alpha_{33}&0\\
\alpha_{41}&\alpha_{42}&\alpha_{43}&0\\
\end{array}
\right)
,\,\left(\begin{array}{ccccc}
0&0&0&0\\
0&0&0&0\\
0&0&0&0\\
\alpha_{41}&\alpha_{42}&\alpha_{43}&\alpha_{44}\\
\end{array}
\right)
,\,
\left(\begin{array}{ccccc}
0&0&0&0\\
0&\alpha_{22}&0&0\\
\alpha_{31}&\alpha_{32}&0&0\\
\alpha_{41}&\alpha_{42}&0&0\\
\end{array}
\right)
,\,
\left(\begin{array}{ccccc}
0&0&0&0\\
\alpha_{21}&\alpha_{22}&0&0\\
\alpha_{31}&\alpha_{32}&0&0\\
\alpha_{41}&\alpha_{42}&0&0\\
\end{array}
\right)
$

\item[$L_8$ : ]
    $\left(\begin{array}{ccccc}
0&0&0&0\\
\alpha_{21}&\alpha_{22}&0&0\\
\alpha_{31}&\alpha_{32}&0&0\\
\alpha_{41}&\alpha_{42}&0&0\\
\end{array}
\right),\,
\left(\begin{array}{ccccc}
0&0&0&0\\
0&\alpha_{22}&0&0\\
\alpha_{31}&\alpha_{32}&0&0\\
\alpha_{41}&\alpha_{42}&0&0\\
\end{array}
\right)
,\,
\left(\begin{array}{ccccc}
0&0&0&0\\
0&0&0&0\\
\alpha_{31}&\alpha_{32}&\alpha_{33}&0\\
\alpha_{41}&\alpha_{42}&\alpha_{43}&0\\
\end{array}
\right)
,\,
\left(\begin{array}{ccccc}
\alpha_{11}&0&0&0\\
0&\alpha_{11}&0&0\\
0&0&\alpha_{11}&0\\
\alpha_{41}&\alpha_{42}&0&\alpha_{11}\\
\end{array}
\right)
,
\\
\left(\begin{array}{ccccc}
0&0&0&0\\
0&0&0&0\\
0&0&0&0\\
\alpha_{41}&\alpha_{42}&\alpha_{43}&\alpha_{44}\\
\end{array}
\right)
$

\item[$L_{9}$ : ]
    $\left(\begin{array}{ccccc}
0&0&0&0\\
\alpha_{21}&\alpha_{22}&0&0\\
\alpha_{31}&\alpha_{32}&0&0\\
\alpha_{41}&\alpha_{42}&0&0\\
\end{array}
\right),\,
\left(\begin{array}{ccccc}
0&0&0&0\\
\alpha_{21}&\alpha_{22}&0&0\\
\alpha_{31}&\alpha_{32}&0&0\\
\alpha_{41}&\alpha_{42}&\alpha_{21}-\alpha_{31}&\alpha_{22}\\
\end{array}
\right)
,\,
\left(\begin{array}{ccccc}
0&0&0&0\\
0&0&0&0\\
0&0&0&0\\
\alpha_{41}&\alpha_{42}&\alpha_{43}&\alpha_{44}\\
\end{array}
\right)
,
\\
\left(\begin{array}{ccccc}
\alpha_{11}&0&0&0\\
0&\alpha_{11}&0&0\\
0&0&\alpha_{11}&0\\
\alpha_{41}&\alpha_{42}&0&\alpha_{11}\\
\end{array}
\right)
;
$

\item[$L_{10}$ : ]
    $\left(\begin{array}{ccccc}
0&0&0&0\\
\alpha_{21}&\alpha_{22}&0&0\\
\alpha_{31}&\alpha_{32}&0&0\\
\alpha_{41}&\alpha_{42}&0&0\\
\end{array}
\right),\,
\left(\begin{array}{ccccc}
0&0&0&0\\
0&0&0&0\\
\alpha_{31}&\alpha_{32}&0&0\\
\alpha_{41}&\alpha_{42}&\alpha_{43}&0\\
\end{array}
\right)
,\,
\left(\begin{array}{ccccc}
0&0&0&0\\
0&\alpha_{22}&0&0\\
0&0&0&0\\
\alpha_{41}&\alpha_{42}&0&\alpha_{22}\\
\end{array}
\right)
,\,
\left(\begin{array}{ccccc}
\alpha_{11}&0&0&0\\
0&\alpha_{11}&0&0\\
0&0&\alpha_{11}&0\\
\alpha_{41}&\alpha_{42}&0&\alpha_{11}\\
\end{array}
\right)
;
$

\item[$L_{11}$ : ]
    $\left(\begin{array}{ccccc}
\alpha_{11}&\alpha_{12}&0&0\\
0&0&0&0\\
\alpha_{31}&\alpha_{32}&0&0\\
\alpha_{41}&\alpha_{42}&\alpha_{43}&\alpha_{11}\\
\end{array}
\right),\,
\left(\begin{array}{ccccc}
0&0&0&0\\
0&\alpha_{22}&0&0\\
\alpha_{31}&\alpha_{32}&0&-2\alpha_{22}\\
\alpha_{41}&\alpha_{42}&0&\alpha_{22}\\
\end{array}
\right)
,\,
\left(\begin{array}{ccccc}
0&0&0&0\\
0&0&0&0\\
\alpha_{31}&\alpha_{32}&\alpha_{33}&\alpha_{34}\\
\alpha_{41}&\alpha_{42}&\alpha_{43}&\alpha_{44}\\
\end{array}
\right),
$

$
\left(\begin{array}{ccccc}
\alpha_{22}&0&0&0\\
0&\alpha_{22}&0&0\\
\alpha_{31}&\alpha_{32}&\alpha_{22}&0\\
\alpha_{41}&\alpha_{42}&0&\alpha_{22}\\
\end{array}
\right),\quad
\left(\begin{array}{ccccc}
2\alpha_{21}+\alpha_{22}&\alpha_{12}&0&0\\
\alpha_{21}&\alpha_{22}&0&0\\
\alpha_{31}&\alpha_{32}&\alpha_{22}&\alpha_{34}\\
\alpha_{41}&\alpha_{42}&\alpha_{43}&2\alpha_{21}+\alpha_{22}
\end{array}
\right)
;
$

\item[$L_{12}$ : ]
    $\left(\begin{array}{ccccc}
0&0&0&0\\
0&0&0&0\\
\alpha_{31}&\alpha_{32}&\alpha_{33}&0\\
\alpha_{41}&\alpha_{42}&0&\alpha_{44}\\
\end{array}
\right),\,
\left(\begin{array}{ccccc}
0&0&0&0\\
0&0&0&0\\
\alpha_{31}&\alpha_{32}&\alpha_{33}&\alpha_{34}\\
\alpha_{41}&\alpha_{42}&\alpha_{43}&\alpha_{44}\\
\end{array}
\right)
,\,
\left(\begin{array}{ccccc}
0&0&0&0\\
0&0&0&0\\
\alpha_{31}&\alpha_{32}&0&\alpha_{34}\\
\alpha_{41}&\alpha_{42}&\alpha_{43}&\alpha_{44}\\
\end{array}
\right),
\,
\left(\begin{array}{ccccc}
\alpha_{22}&0&0&0\\
0&\alpha_{22}&0&0\\
\alpha_{31}&\alpha_{32}&\alpha_{22}&0\\
\alpha_{41}&\alpha_{42}&0&\alpha_{22}\\
\end{array}
\right)
;
$

\item[$L_{13}$ : ]
    $\left(\begin{array}{ccccc}
0&0&0&0\\
0&0&0&0\\
\alpha_{31}&\alpha_{32}&\alpha_{33}&0\\
\alpha_{41}&\alpha_{42}&\alpha_{43}&\alpha_{44}\\
\end{array}
\right),\,
\left(\begin{array}{ccccc}
0&0&0&0\\
0&0&0&0\\
\alpha_{31}&\alpha_{32}&0&\alpha_{34}\\
\alpha_{41}&\alpha_{42}&\alpha_{43}&\alpha_{44}\\
\end{array}
\right)
,\,
\left(\begin{array}{ccccc}
0&0&0&0\\
0&0&0&0\\
\alpha_{31}&\alpha_{32}&0&\alpha_{34}\\
\alpha_{41}&\alpha_{42}&\alpha_{43}&0\\
\end{array}
\right),
\left(\begin{array}{ccccc}
\alpha_{22}&0&0&0\\
0&\alpha_{22}&0&0\\
\alpha_{31}&\alpha_{32}&\alpha_{22}&0\\
\alpha_{41}&\alpha_{42}&0&\alpha_{22}\\
\end{array}
\right)
;
$

\item[$L_{14}$ : ]
    $\left(\begin{array}{ccccc}
0&0&0&0\\
0&0&0&0\\
0&0&\alpha_{44}&0\\
\alpha_{41}&\alpha_{42}&\alpha_{43}&\alpha_{44}\\
\end{array}
\right),
\left(\begin{array}{ccccc}
\alpha_{44}&0&0&0\\
0&\alpha_{44}&0&0\\
0&0&0&0\\
\alpha_{41}&\alpha_{42}&\alpha_{43}&\alpha_{44}\\
\end{array}
\right),
\left(\begin{array}{ccccc}
\alpha_{44}&0&0&0\\
0&\alpha_{44}&0&0\\
0&0&\alpha_{44}&0\\
\alpha_{41}&\alpha_{42}&0&\alpha_{44}\\
\end{array}
\right),
\left(\begin{array}{ccccc}
\alpha_{44}&0&0&0\\
0&\alpha_{44}&0&0\\
0&0&\alpha_{33}&0\\
\alpha_{41}&\alpha_{42}&0&\alpha_{44}\\
\end{array}
\right)
;
$

\item[$L_{15}$ : ]
    $\left(\begin{array}{ccccc}
0&0&0&0\\
0&0&0&0\\
\alpha_{31}&\alpha_{32}&\alpha_{33}&0\\
\alpha_{41}&\alpha_{42}&\alpha_{43}&0\\
\end{array}
\right),
\left(\begin{array}{ccccc}
\alpha_{11}&\alpha_{12}&\alpha_{13}&0\\
0&0&0&0\\
0&0&0&0\\
\alpha_{41}&\alpha_{42}&\alpha_{43}&0\\
\end{array}
\right),
\left(\begin{array}{ccccc}
0&0&0&0\\
0&0&0&0\\
0&0&0&0\\
\alpha_{41}&\alpha_{42}&\alpha_{43}&\alpha_{44}\\
\end{array}
\right),
\left(\begin{array}{ccccc}
\alpha_{44}&0&0&0\\
0&\alpha_{44}&0&0\\
0&0&\alpha_{44}&0\\
\alpha_{41}&\alpha_{42}&\alpha_{43}&\alpha_{44}\\
\end{array}
\right)
;
$

\item[$L_{16}$ : ]
    $\left(\begin{array}{ccccc}
0&0&0&0\\
\alpha_{21}&\alpha_{22}&\alpha_{23}&0\\
0&0&0&0\\
\alpha_{41}&\alpha_{42}&\alpha_{43}&0\\
\end{array}
\right),
\left(\begin{array}{ccccc}
\alpha_{11}&0&0&0\\
0&\alpha_{11}&0&0\\
0&0&0&0\\
\alpha_{41}&\alpha_{42}&\alpha_{43}&\alpha_{11}\\
\end{array}
\right),
\left(\begin{array}{ccccc}
0&0&0&0\\
0&0&0&0\\
0&0&\alpha_{44}&0\\
\alpha_{41}&\alpha_{42}&\alpha_{43}&\alpha_{44}\\
\end{array}
\right),
\left(\begin{array}{ccccc}
\alpha_{11}&0&0&0\\
0&\alpha_{11}&0&0\\
0&0&\alpha_{11}&0\\
\alpha_{41}&\alpha_{42}&\alpha_{43}&\alpha_{11}\\
\end{array}
\right)
;
$

\item[$L_{17}$ : ]
    $\left(\begin{array}{ccccc}
0&0&0&0\\
0&0&0&0\\
\alpha_{31}&\alpha_{32}&\alpha_{33}&\alpha_{34}\\
\alpha_{41}&\alpha_{42}&\alpha_{43}&\alpha_{44}\\
\end{array}
\right),
\left(\begin{array}{ccccc}
\alpha_{11}&\alpha_{12}&0&0\\
0&0&0&0\\
\alpha_{31}&\alpha_{32}&0&0\\
\alpha_{41}&\alpha_{42}&\alpha_{43}&\alpha_{44}\\
\end{array}
\right),
\left(\begin{array}{ccccc}
0&0&0&0\\
0&0&0&0\\
\alpha_{31}&\alpha_{32}&0&0\\
\alpha_{41}&\alpha_{42}&\alpha_{43}&0\\
\end{array}
\right),
\left(\begin{array}{ccccc}
\alpha_{11}&0&0&0\\
0&\alpha_{11}&0&0\\
\alpha_{31}&\alpha_{32}&\alpha_{11}&0\\
\alpha_{41}&\alpha_{42}&\alpha_{43}&\alpha_{11}\\
\end{array}
\right)
;
$

\item[$L_{18}$ : ]
    $\left(\begin{array}{ccccc}
0&0&0&0\\
0&0&0&0\\
\alpha_{31}&\alpha_{32}&\alpha_{33}&\alpha_{34}\\
\alpha_{41}&\alpha_{42}&\alpha_{43}&\alpha_{44}\\
\end{array}
\right),
\left(\begin{array}{ccccc}
\alpha_{11}&0&0&0\\
0&0&0&0\\
\alpha_{31}&\alpha_{32}&0&0\\
\alpha_{41}&\alpha_{42}&0&0\\
\end{array}
\right),
\left(\begin{array}{ccccc}
\alpha_{11}&\alpha_{12}&0&0\\
0&0&0&0\\
\alpha_{31}&\alpha_{32}&0&0\\
\alpha_{41}&\alpha_{42}&0&0\\
\end{array}
\right),
\left(\begin{array}{ccccc}
\alpha_{22}&0&0&0\\
0&\alpha_{22}&0&0\\
\alpha_{31}&\alpha_{32}&\alpha_{22}&0\\
\alpha_{41}&\alpha_{42}&0&\alpha_{22}\\
\end{array}
\right)
;
$

\item[$L_{19}$ : ]
    $\left(\begin{array}{ccccc}
0&0&0&0\\
0&0&0&0\\
\alpha_{31}&\alpha_{32}&\alpha_{33}&\alpha_{34}\\
\alpha_{41}&\alpha_{42}&\alpha_{43}&\alpha_{44}\\
\end{array}
\right),
\left(\begin{array}{ccccc}
\alpha_{11}&\alpha_{12}&0&0\\
0&0&0&0\\
\alpha_{31}&\alpha_{32}&\alpha_{33}&0\\
\alpha_{41}&\alpha_{42}&\alpha_{43}&0\\
\end{array}
\right),
\left(\begin{array}{ccccc}
0&\alpha_{12}&\alpha_{13}&0\\
0&0&0&0\\
\alpha_{31}&\alpha_{32}&\alpha_{33}&0\\
\alpha_{41}&\alpha_{42}&\alpha_{43}&0\\
\end{array}
\right),
\left(\begin{array}{ccccc}
\alpha_{22}&0&0&0\\
0&\alpha_{22}&0&0\\
\alpha_{31}&\alpha_{32}&\alpha_{22}&0\\
\alpha_{41}&\alpha_{42}&0&\alpha_{22}\\
\end{array}
\right)
;
$

\item[$L_{20}$ : ]
    $\left(\begin{array}{ccccc}
0&0&0&0\\
0&0&0&0\\
\alpha_{31}&\alpha_{32}&\alpha_{33}&\alpha_{34}\\
\alpha_{41}&\alpha_{42}&\alpha_{43}&\alpha_{44}\\
\end{array}
\right),
\left(\begin{array}{ccccc}
\alpha_{11}&\alpha_{12}&\alpha_{13}&0\\
0&0&0&0\\
\alpha_{31}&\alpha_{32}&\alpha_{33}&0\\
\alpha_{41}&\alpha_{42}&\alpha_{43}&0\\
\end{array}
\right),
\left(\begin{array}{ccccc}
0&0&0&0\\
0&\alpha_{22}&0&0\\
\alpha_{31}&\alpha_{32}&\alpha_{22}&0\\
\alpha_{41}&\alpha_{42}&0&0\\
\end{array}
\right),
\left(\begin{array}{ccccc}
\alpha_{22}&0&0&0\\
0&\alpha_{22}&0&0\\
\alpha_{31}&\alpha_{32}&\alpha_{22}&0\\
\alpha_{41}&\alpha_{42}&0&\alpha_{22}\\
\end{array}
\right)
;
$

\item[$L_{21}$ : ]
    $\left(\begin{array}{ccccc}
0&0&0&0\\
\alpha_{21}&\alpha_{22}&\alpha_{23}&0\\
0&0&0&0\\
\alpha_{41}&\alpha_{42}&\alpha_{43}&0\\
\end{array}
\right),\quad
\left(\begin{array}{ccccc}
\alpha_{44}&0&0&0\\
0&\alpha_{44}&0&0\\
0&0&0&0\\
\alpha_{41}&\alpha_{42}&\alpha_{43}&\alpha_{44}\\
\end{array}
\right),\quad
\left(\begin{array}{ccccc}
\alpha_{44}&0&0&0\\
0&\alpha_{44}&0&0\\
0&0&\alpha_{44}&0\\
\alpha_{41}&\alpha_{42}&\alpha_{43}&\alpha_{44}\\
\end{array}
\right)
;
$
\end{itemize}
\end{prop}

\subsection{Averaging operator}\,

Let $(L, \llbracket-,-\rrbracket)$   be an $n$-dimensional Leibniz algebra with basis $\{e_i\}\, (1 \leq i \leq n)$ and let $\beta$ be the averaging operator on $L$.
For any $i,j,k\in\mathbb{N}$, $1 \leq i,j,k \leq n$, let us put
$$\llbracket e_i,e_j\rrbracket=\sum_{k=1}^nc_{ij}^ke_k\quad \text{and}\quad \beta(e_i)=\sum_{j=1}^n\beta_{ji}e_j.$$
Then, in term of basis elements, equation (\ref{eq6}) is equivalent to

\begin{align*} 
   \sum^{n}_{p=1} \sum^{n}_{k=1}\beta_{ki}C_{kj}^{p}\beta_{qp}=\sum^{n}_{k=1}\sum^{n}_{p=1}\beta_{ki}\beta_{pj} C_{kp}^{p}=
   \sum^{n}_{p=1} \sum^{n}_{k=1}\beta_{kj}C_{ik}^{p}\beta_{qp},
\end{align*} for $i,j,q=1,2,\dots,n$.

\begin{prop}
The description of the Averaging operator of every 4-dimensional  Leibniz algebra is given below.
\begin{itemize}
    \item[$L_1$ : ]
    $\left(\begin{array}{ccccc}
\beta_{11}&0&0&0\\
0&\beta_{11}&0&0\\
0&0&\beta_{11}&0\\
\beta_{41}&0&0&\beta_{11}\\
\end{array}
\right),\,
\left(\begin{array}{ccccc}
0&0&0&0\\
0&0&0&0\\
\beta_{31}&\beta_{32}&\beta_{33}&0\\
\beta_{41}&\beta_{42}&\beta_{43}&0\\
\end{array}
\right)
,\,
\left(\begin{array}{ccccc}
0&0&0&0\\
\beta_{21}&\beta_{22}&0&\beta_{24}\\
0&0&0&0\\
\beta_{41}&\beta_{42}&0&\beta_{44}\\
\end{array}
\right)
,\,
\left(\begin{array}{ccccc}
0&0&0&0\\
0&0&0&0\\
\beta_{21}&\beta_{32}&\beta_{33}&0\\
\beta_{41}&\beta_{42}&\beta_{43}&0\\
\end{array}
\right)
;
$

 \item[$L_2$ : ]
    $\left(\begin{array}{ccccc}
\beta_{11}&0&0&0\\
0&\beta_{11}&0&0\\
0&0&\beta_{11}&0\\
\beta_{41}&\beta_{42}&0&\beta_{11}\\
\end{array}
\right),\,
\left(\begin{array}{ccccc}
0&0&0&0\\
0&0&0&0\\
0&0&0&0\\
\beta_{41}&\beta_{42}&\beta_{43}&\beta_{44}\\
\end{array}
\right)
,\,
\left(\begin{array}{ccccc}
0&0&0&0\\
0&0&0&0\\
\beta_{31}&\beta_{32}&\beta_{33}&0\\
\beta_{41}&\beta_{42}&\beta_{43}&0\\
\end{array}
\right)
,\,
\left(\begin{array}{ccccc}
0&0&0&0\\
0&\beta_{22}&0&0\\
\beta_{31}&\beta_{32}&0&0\\
\beta_{41}&\beta_{42}&0&0\\
\end{array}
\right)
;
$

\item[$L_3$ : ]
    $\left(\begin{array}{ccccc}
\beta_{11}&0&0&0\\
0&\beta_{11}&0&0\\
0&0&\beta_{11}&0\\
\beta_{41}&\beta_{42}&0&\beta_{11}\\
\end{array}
\right),\,
\left(\begin{array}{ccccc}
0&0&0&0\\
\beta_{21}&\beta_{22}&0&0\\
0&0&0&0\\
\beta_{41}&\beta_{42}&0&\beta_{44}\\
\end{array}
\right)
,\,
\left(\begin{array}{ccccc}
0&0&0&0\\
\beta_{21}&\beta_{22}&\beta_{24}&0\\
0&0&0&0\\
\beta_{41}&\beta_{42}&0&\beta_{44}\\
\end{array}
\right)
,\,
\left(\begin{array}{ccccc}
0&0&0&0\\
0&\beta_{22}&0&0\\
0&0&0&0\\
\beta_{41}&\beta_{42}&0&\beta_{44}\\
\end{array}
\right)
;
$

\item[$L_4$ : ]
    $\left(\begin{array}{ccccc}
\beta_{11}&0&0&0\\
0&\beta_{11}&0&0\\
0&0&\beta_{11}&0\\
\beta_{41}&\beta_{42}&0&\beta_{11}\\
\end{array}
\right)
,\,
\left(\begin{array}{ccccc}
0&0&0&0\\
0&0&0&0\\
\beta_{31}&\beta_{32}&\beta_{33}&0\\
\beta_{41}&\beta_{42}&\beta_{43}&0\\
\end{array}
\right)
,\,
\left(\begin{array}{ccccc}
0&0&0&0\\
0&0&0&0\\
0&0&0&0\\
\beta_{41}&\beta_{42}&\beta_{43}&\beta_{44}\\
\end{array}
\right)
;
$

\item[$L_5$ : ]
    $\left(\begin{array}{ccccc}
\beta_{11}&0&0&0\\
0&\beta_{11}&0&0\\
0&0&\beta_{11}&0\\
\beta_{41}&\beta_{42}&0&\beta_{11}\\
\end{array}
\right),
\left(\begin{array}{ccccc}
0&0&0&0\\
0&0&0&0\\
\beta_{31}&\beta_{32}&\beta_{33}\\
\beta_{41}&\beta_{42}&\beta_{43}\\
\end{array}
\right),
\left(\begin{array}{ccccc}
\beta_{11}&0&0&0\\
-\beta_{11}&0&0&0\\
-\beta_{11}&-\beta_{11}&-\beta_{11}&0\\
\beta_{41}&\beta_{42}&0&\beta_{11}\\
\end{array}
\right),
\left(\begin{array}{ccccc}
0&0&0&0\\
0&0&0&0\\
0&0&0&0\\
\beta_{41}&\beta_{42}&\beta_{43}&\beta_{44}\\
\end{array}
\right)
;
$

\item[$L_6$ : ]
    $\left(\begin{array}{ccccc}
\beta_{11}&0&0&0\\
0&\beta_{11}&0&0\\
0&0&\beta_{11}&0\\
\beta_{41}&\beta_{42}&0&\beta_{11}\\
\end{array}
\right),
\left(\begin{array}{ccccc}
\beta_{11}&0&0&0\\
0&0&0&0\\
0&0&\beta_{11}&0\\
\beta_{41}&\beta_{42}&0&\beta_{11}\\
\end{array}
\right),
\left(\begin{array}{ccccc}
0&0&0&0\\
\beta_{41}&\beta_{42}&\beta_{43}&\beta_{22}\\
0&0&0&0\\
\beta_{41}&\beta_{42}&0&\beta_{11}\\
\end{array}
\right),
\left(\begin{array}{ccccc}
0&0&0&0\\
0&0&0&0\\
0&0&0&0\\
\beta_{41}&\beta_{42}&\beta_{43}&\beta_{44}\\
\end{array}
\right)
;
$

\item[$L_7$ : ]
    $\left(\begin{array}{ccccc}
\beta_{11}&0&0&0\\
0&\beta_{11}&0&0\\
0&0&\beta_{11}&0\\
\beta_{41}&\beta_{42}&0&\beta_{11}\\
\end{array}
\right),
\left(\begin{array}{ccccc}
0&0&0&0\\
0&0&\beta_{22}&0\\
0&0&0&0\\
\beta_{41}&\beta_{42}&0&\beta_{44}\\
\end{array}
\right),
\left(\begin{array}{ccccc}
0&0&0&0\\
0&0&0&0\\
\beta_{31}&\beta_{32}&\beta_{33}&0\\
\beta_{41}&\beta_{42}&\beta_{43}&0\\
\end{array}
\right),
\left(\begin{array}{ccccc}
0&0&0&0\\
\beta_{21}&\beta_{22}&0&0\\
\beta_{31}&\beta_{32}&0&0\\
\beta_{41}&\beta_{42}&0&0\\
\end{array}
\right)
;
$

\item[$L_8$ : ]
    $\left(\begin{array}{ccccc}
\beta_{11}&0&0&0\\
0&\beta_{11}&0&0\\
0&0&\beta_{11}&0\\
\beta_{41}&\beta_{42}&0&\beta_{11}\\
\end{array}
\right),
\left(\begin{array}{ccccc}
0&0&0&0\\
0&0&\beta_{22}&0\\
0&0&0&0\\
\beta_{41}&\beta_{42}&\beta_{43}&\beta_{44}\\
\end{array}
\right),
\left(\begin{array}{ccccc}
0&0&0&0\\
0&0&0&0\\
\beta_{31}&\beta_{32}&\beta_{33}&0\\
\beta_{41}&\beta_{42}&\beta_{43}&0\\
\end{array}
\right),
\left(\begin{array}{ccccc}
0&0&0&0\\
\beta_{21}&0&0&0\\
\beta_{31}&\beta_{32}&0&0\\
\beta_{41}&\beta_{42}&0&0\\
\end{array}
\right)
;
$

\item[$L_9$ : ]
    $\left(\begin{array}{ccccc}
\beta_{11}&0&0&0\\
0&\beta_{11}&0&0\\
0&0&\beta_{11}&0\\
\beta_{41}&\beta_{42}&0&\beta_{11}\\
\end{array}
\right),
\left(\begin{array}{ccccc}
0&0&0&0\\
\beta_{21}&\beta_{21}&0&0\\
\beta_{31}&\beta_{32}&0&0\\
\beta_{41}&\beta_{42}&\beta_{21}-\beta_{31}&\beta_{44}\\
\end{array}
\right),
\left(\begin{array}{ccccc}
0&0&0&0\\
0&0&0&0\\
\beta_{31}&\beta_{32}&\beta_{33}&0\\
\beta_{41}&\beta_{42}&\beta_{43}&0\\
\end{array}
\right),
\left(\begin{array}{ccccc}
0&0&0&0\\
\beta_{21}&0&0&0\\
\beta_{31}&\beta_{32}&0&0\\
\beta_{41}&\beta_{42}&0&0\\
\end{array}
\right)
;
$

\item[$L_{10}$ : ]
    $\left(\begin{array}{ccccc}
\beta_{11}&0&0&0\\
0&\beta_{11}&0&0\\
0&0&\beta_{11}&0\\
\beta_{41}&\beta_{42}&0&\beta_{11}\\
\end{array}
\right),
\left(\begin{array}{ccccc}
0&0&0&0\\
0&0&0&0\\
\beta_{31}&\beta_{32}&0&0\\
\beta_{41}&\beta_{42}&\beta_{43}&0\
\end{array}
\right),
\left(\begin{array}{ccccc}
0&0&0&0\\
0&0&0&0\\
0&0&\beta_{33}&0\\
\beta_{41}&\beta_{42}&0&\beta_{22}\\
\end{array}
\right),
\left(\begin{array}{ccccc}
0&0&0&0\\
\beta_{21}&\beta_{22}&0&0\\
\frac{\beta_{21}^2}{\beta_{22}}&\beta_{31}&0&0\\
\beta_{41}&\beta_{42}&0&\beta_{22}\\
\end{array}
\right)
;
$

\item[$L_{11}$ : ]
    $\left(\begin{array}{ccccc}
\beta_{22}&0&0&0\\
0&\beta_{22}&0&0\\
\beta_{31}&\beta_{32}&\beta_{22}&0\\
\beta_{41}&\beta_{42}&0&\beta_{22}\\
\end{array}
\right),
\left(\begin{array}{ccccc}
\beta_{11}&0&0&0\\
0&0&0&0\\
\beta_{31}&\beta_{32}&0&0\\
\beta_{41}&\beta_{42}&0&\beta_{11}\\
\end{array}
\right),
\left(\begin{array}{ccccc}
0&0&0&0\\
0&0&0&0\\
\beta_{31}&\beta_{32}&0&0\\
\beta_{41}&\beta_{42}&0&\beta_{11}\\
\end{array}
\right),
\left(\begin{array}{ccccc}
0&0&0&0\\
0&\beta_{22}&0&0\\
\beta_{31}&\beta_{32}&0&-2\beta_{22}\\
\beta_{41}&\beta_{42}&0&\beta_{22}\\
\end{array}
\right)
;
$

\item[$L_{12}$ : ]
    $\left(\begin{array}{ccccc}
\beta_{22}&0&0&0\\
0&\beta_{22}&0&0\\
\beta_{31}&\beta_{32}&\beta_{22}&\beta_{34}\\
\beta_{41}&\beta_{42}&0&\beta_{44}\\
\end{array}
\right),
\left(\begin{array}{ccccc}
0&\beta_{12}&0&0\\
0&0&0&0\\
\beta_{31}&0&\beta_{34}&0\\
\beta_{41}&\beta_{42}&0&\beta_{44}\\
\end{array}
\right),
\left(\begin{array}{ccccc}
0&0&0&0\\
0&0&0&0\\
\beta_{31}&\beta_{32}&\beta_{33}&\beta_{34}\\
\beta_{41}&\beta_{42}&0&\beta_{44}\\
\end{array}
\right),
\left(\begin{array}{ccccc}
0&0&0&0\\
0&0&0&0\\
\beta_{31}&\beta_{32}&\beta_{33}&\beta_{34}\\
\beta_{41}&\beta_{42}&0&\beta_{44}\\
\end{array}
\right)
;
$

\item[$L_{13}$ : ]
    $\left(\begin{array}{ccccc}
\beta_{11}&0&0&0\\
0&\beta_{11}&0&0\\
\beta_{31}&\beta_{32}&\beta_{11}&0\\
\beta_{41}&\beta_{42}&0&\beta_{11}\\
\end{array}
\right),
\left(\begin{array}{ccccc}
0&0&0&0\\
0&0&0&0\\
\beta_{31}&\beta_{32}&\beta_{33}&\beta_{34}\\
\beta_{41}&\beta_{42}&\beta_{43}&\beta_{44}\\
\end{array}
\right),
\left(\begin{array}{ccccc}
0&0&0&0\\
\beta_{21}&\frac{\beta_{21}}{g}&0&0\\
\beta_{31}&\beta_{32}&\frac{\beta_{21}}{g}&0\\
\beta_{41}&\beta_{42}&0&-\beta_{21}\\
\end{array}
\right)
;
$

\item[$L_{14}$ : ]
    $\left(\begin{array}{ccccc}
\beta_{44}&0&0&0\\
0&\beta_{44}&0&0\\
0&0&\beta_{33}&0\\
\beta_{41}&\beta_{42}&\beta_{43}&\beta_{44}\\
\end{array}
\right),
\left(\begin{array}{ccccc}
0&0&0&0\\
0&0&0&0\\
0&0&0&0\\
\beta_{41}&\beta_{42}&\beta_{43}&\beta_{44}\\
\end{array}
\right),
\left(\begin{array}{ccccc}
0&0&0&0\\
0&0&0&0\\
0&0&\beta_{44}&0\\
\beta_{41}&\beta_{42}&\beta_{43}&\beta_{44}\\
\end{array}
\right)
,
\left(\begin{array}{ccccc}
\beta_{44}&0&0&0\\
0&\beta_{44}&0&0\\
0&0&\beta_{44}&0\\
\beta_{41}&\beta_{42}&\beta_{43}&\beta_{44}\\
\end{array}
\right)
;
$

\item[$L_{15}$ : ]
    $\left(\begin{array}{ccccc}
\beta_{22}&0&0&0\\
0&\beta_{22}&0&0\\
0&0&\beta_{22}&0\\
\beta_{41}&\beta_{42}&0&\beta_{22}\\
\end{array}
\right),
\left(\begin{array}{ccccc}
0&0&0&0\\
0&0&0&0\\
0&0&0&0\\
\beta_{41}&\beta_{42}&\beta_{43}&\beta_{44}\\
\end{array}
\right),
\left(\begin{array}{ccccc}
0&0&0&0\\
0&0&0&0\\
\beta_{31}&\beta_{32}&0&0\\
\beta_{41}&\beta_{42}&0&0\\
\end{array}
\right)
,
\left(\begin{array}{ccccc}
0&0&0&0\\
0&\beta_{22}&0&0\\
0&0&0&0\\
\beta_{41}&\beta_{42}&-\beta_{22}&\beta_{22}\\
\end{array}
\right)
;
$

\item[$L_{16}$ : ]
    $\left(\begin{array}{ccccc}
\beta_{11}&0&0&0\\
0&\beta_{11}&0&0\\
0&0&\beta_{11}&0\\
\beta_{41}&\beta_{42}&\beta_{43}&\beta_{11}\\
\end{array}
\right),
\left(\begin{array}{ccccc}
0&0&0&0\\
0&0&0&0\\
0&0&0&0\\
\beta_{41}&\beta_{42}&\beta_{43}&\beta_{44}\\
\end{array}
\right),
\left(\begin{array}{ccccc}
\beta_{11}&0&0&0\\
0&\beta_{11}&0&0\\
0&0&0&0\\
\beta_{41}&\beta_{42}&\beta_{43}&\beta_{11}\\
\end{array}
\right)
,
\left(\begin{array}{ccccc}
0&0&0&0\\
\beta_{21}&\beta_{22}&\beta_{23}&0\\
0&0&0&0\\
\beta_{41}&\beta_{42}&\beta_{43}&0\\
\end{array}
\right)
;
$

\item[$L_{17}$ : ]
    $\left(\begin{array}{ccccc}
\beta_{11}&0&0&0\\
0&\beta_{11}&0&0\\
\beta_{31}&\beta_{32}&\beta_{11}&0\\
\beta_{41}&\beta_{42}&0&\beta_{11}\\
\end{array}
\right),
\left(\begin{array}{ccccc}
0&0&0&0\\
0&0&0&0\\
\beta_{31}&\beta_{32}&\beta_{33}&\beta_{34}\\
\beta_{41}&\beta_{42}&\beta_{43}&\beta_{44}\\
\end{array}
\right),
\left(\begin{array}{ccccc}
\beta_{11}&0&0&0\\
0&0&0&0\\
\beta_{31}&\beta_{32}&0&0\\
\beta_{41}&\beta_{42}&0&0\\
\end{array}
\right)
,
\left(\begin{array}{ccccc}
0&0&0&0\\
0&0&0&0\\
\beta_{31}&\beta_{32}&0&\beta_{34}\\
\beta_{41}&\beta_{42}&\beta_{43}&0\\
\end{array}
\right)
;
$

\item[$L_{18}$ : ]
    $\left(\begin{array}{ccccc}
\beta_{22}&0&0&0\\
0&\beta_{22}&0&0\\
\beta_{31}&\beta_{32}&\beta_{22}&0\\
\beta_{41}&\beta_{42}&0&\beta_{22}\\
\end{array}
\right),
\left(\begin{array}{ccccc}
\beta_{11}&\beta_{12}&0&\beta_{14}\\
0&0&0&0\\
\beta_{31}&\beta_{32}&0&\beta_{34}\\
\beta_{41}&\beta_{42}&0&\beta_{44}\\
\end{array}
\right),
\left(\begin{array}{ccccc}
0&0&0&0\\
0&0&0&0\\
\beta_{31}&\beta_{32}&\beta_{33}&\beta_{34}\\
\beta_{41}&\beta_{42}&\beta_{43}&\beta_{44}\\
\end{array}
\right)
,
\left(\begin{array}{ccccc}
0&0&0&0\\
0&0&0&0\\
\beta_{31}&\beta_{32}&\beta_{33}&\beta_{34}\\
\beta_{41}&\beta_{42}&0&\beta_{44}\\
\end{array}
\right)
;
$

\item[$L_{19}$ : ]
    $\left(\begin{array}{ccccc}
\beta_{22}&0&0&0\\
0&\beta_{22}&0&0\\
\beta_{31}&\beta_{32}&\beta_{22}&0\\
\beta_{41}&\beta_{42}&0&\beta_{22}\\
\end{array}
\right),
\left(\begin{array}{ccccc}
\beta_{11}&\beta_{12}&0&0\\
0&0&0&0\\
\beta_{31}&\beta_{32}&\beta_{33}&0\\
\beta_{41}&\beta_{42}&\beta_{43}&0\\
\end{array}
\right),
\left(\begin{array}{ccccc}
0&\beta_{12}&0&0\\
0&\beta_{22}&0&0\\
\beta_{31}&\beta_{32}&\beta_{22}&0\\
\beta_{41}&\beta_{42}&\beta_{43}&0\\
\end{array}
\right)
,
\left(\begin{array}{ccccc}
0&0&0&0\\
0&0&0&0\\
0&0&0&0\\
\beta_{41}&\beta_{42}&\beta_{43}&\beta_{44}\\
\end{array}
\right)
;
$

\item[$L_{20}$ : ]
    $\left(\begin{array}{ccccc}
\beta_{22}&0&0&0\\
0&\beta_{22}&0&0\\
\beta_{31}&\beta_{32}&\beta_{33}&0\\
\beta_{41}&\beta_{42}&\beta_{43}&\beta_{22}\\
\end{array}
\right),
\left(\begin{array}{ccccc}
\beta_{11}&0&\beta_{13}&0\\
0&\beta_{11}&0&0\\
\beta_{31}&\beta_{32}&\beta_{33}&0\\
\beta_{41}&\beta_{42}&\beta_{43}&\beta_{11}\\
\end{array}
\right),
\left(\begin{array}{ccccc}
0&0&0&0\\
0&0&0&0\\
\beta_{31}&\beta_{32}&\beta_{33}&\beta_{34}\\
\beta_{41}&\beta_{42}&\beta_{43}&\beta_{44}\\
\end{array}
\right)
,
\left(\begin{array}{ccccc}
0&0&0&0\\
\beta_{21}&\beta_{22}&\beta_{23}&0\\
\beta_{31}&\beta_{32}&\beta_{33}&0\\
\beta_{41}&\beta_{42}&\beta_{43}&0\\
\end{array}
\right)
;
$

\item[$L_{21}$ : ]
    $\left(\begin{array}{ccccc}
\beta_{44}&0&0&0\\
0&\beta_{44}&0&0\\
0&0&\beta_{44}&0\\
\beta_{41}&\beta_{42}&\beta_{43}&\beta_{44}\\
\end{array}
\right),
\left(\begin{array}{ccccc}
\beta_{44}&0&0&0\\
0&\beta_{44}&0&0\\
0&0&0&0\\
\beta_{41}&\beta_{42}&\beta_{43}&\beta_{44}\\
\end{array}
\right),
\left(\begin{array}{ccccc}
0&0&\beta_{13}&0\\
0&0&\beta_{23}&0\\
0&0&0&0\\
\beta_{41}&\beta_{42}&\beta_{43}&0\\
\end{array}
\right)
,
\left(\begin{array}{ccccc}
0&0&0&0\\
\beta_{21}&\beta_{22}&\beta_{23}&0\\
0&0&0&0\\
\beta_{41}&\beta_{42}&\beta_{43}&0\\
\end{array}
\right).
$
\end{itemize}
\end{prop}
\begin{cor}\,
\begin{itemize}
    \item The dimensions of the Rota-Baxter operator of Leibniz algebras of $4$-dimensional range between $3$ and $10$.
    \item The dimensions of the Nenjenhuis operator of Leibniz algebras of $4$-dimensional range between $5$ and $10$.
    \item The dimensions of the Reynolds operator of Leibniz algebras of $4$-dimensional range between $2$ and $9$.
    \item The dimensions of the Averaging operator of Leibniz algebras of $4$-dimensional range between $2$ and $9$.
\end{itemize}
\end{cor}

\subsection{Compatible Leibniz Algebras}\,

 Let $(L, \llbracket-,-\rrbracket_1)$, $(L, \llbracket-,-\rrbracket_2)$   be an $n$-dimensional compatible Leibniz algebra with basis $\{e_i\}\, (1 \leq i \leq n)$
For any $i,j,k\in\mathbb{N}$, $1 \leq i,j,k \leq n$, let us put
$$\llbracket e_i,e_j\rrbracket_1=\sum_{k=1}^nC_{ij}^ke_k\quad \text{and}\quad \llbracket e_i,e_j\rrbracket_2=\sum_{k=1}^n\gamma_{ij}^ke_k.$$
The axioms in Definition \ref{dcl} is equivalent to
\begin{align}
   \sum^{n}_{p=1}\Bigg(C_{jk}^{p}C_{ip}^{q}-C_{ij}^pC_{pk}^{q}+C_{ik}^{p}C_{pj}^q\Bigg)=0,
\end{align}
\begin{align}
 \sum^{n}_{p=1}\Bigg(\gamma_{jk}^{p}\gamma_{ip}^{q}-\gamma_{ij}^p\gamma_{pk}^{q}+\gamma_{ik}^{p}\gamma_{pj}^q\Bigg)=0,
\end{align}

\begin{align}
   \sum^{n}_{p=1}\Bigg(C_{jk}^{p}\gamma_{ip}^{q}+\gamma_{jk}^{p}C_{ip}^{q}-C_{ij}^p\gamma_{pk}^{q}-\gamma_{ij}^pC_{pk}^{q}+C_{ik}^{p}\gamma_{pj}^q+\gamma_{ik}^{p}C_{pj}^q\Bigg)=0.
\end{align}
\newpage
\begin{thm}
The four-dimenssional compatible right Leibniz algebras are given by the pairs 

    \begin{multicols}{5}
\begin{enumerate}[1.]
     \item $(L_1, L_3)$,
    \item $(L_2, L_4)$,
    \item $(L_2, L_5)$,
    \item $(L_2, L_6)$,
    \item $(L_2,L_{13})$,
     \item $(L_3, L_4)$,
      \item $(L_3, L_5)$,
      \item $(L_3, L_6)$,
      \item $(L_3, L_{13})$,
      \item $(L_4, L_{5})$,
      \item $(L_4, L_{6})$,
      \item $(L_4, L_{9})$,
      \item $(L_4, L_{13})$,
      \item $(L_5, L_{7})$,
      \item $(L_5, L_{13})$,
      \item $(L_{6},L_{13})$,
      \item $(L_7, L_{18})$,
      \item $(L_8, L_{9})$,
      \item $(L_8, L_{10})$,
      \item $(L_8, L_{18})$,
      \item $(L_9, L_{10})$,
      \item $(L_9,L_{18})$,
      \item $(L_{10},L_{18})$,
       \item $(L_{11}, L_{12})$,
     \item $(L_{11}, L_{13})$,
      \item $(L_{11}, L_{17})$,
      \item $(L_{11}, L_{18})$,
     \item $(L_{11}, L_{19})$,
    \item $(L_{11}, L_{20})$,
    \item $(L_{12}, L_{23})$,
    \item $(L_{12}, L_{17})$, 
    \item $(L_{12}, L_{18})$,
    \item $(L_{12}, L_{19})$,
      \item $(L_{12}, L_{20})$,
    \item $(L_{13}, L_{17})$,
      \item $(L_{13}, L_{18})$,
    \item $(L_{13},L_{19})$,
    \item $(L_{13}, L_{20})$, 
    \item $(L_{14}, L_{15})$,
    \item $(L_{14}, L_{16})$,
    \item $(L_{14}, L_{21})$, 
   \item $(L_{15}, L_{16})$,
    \item $(L_{15}, L_{21})$,
    \item $(L_{16}, L_{21})$,  
   \item $(L_{17}, L_{18})$,
   \item $(L_{17}, L_{19})$, 
   \item $(L_{17}, L_{20})$, 
   \item $(L_{18}, L_{19})$, 
  \item $(L_{18}, L_{20})$, 
  \item $(L_{19}, L_{21})$.        
\end{enumerate}
\end{multicols}
\end{thm}

\end{document}